\documentclass[a4paper,fleqn]{cas-sc}

\usepackage{siunitx}
\usepackage[numbers]{natbib}
\usepackage{todonotes}

\def\tsc#1{\csdef{#1}{\textsc{\lowercase{#1}}\xspace}}
\tsc{WGM}
\tsc{QE}

\def\dO{\partial \Omega}

\def\Ol{\overline \Omega}

\def\hX{{\hat X}}

\def\scriptO{{{\it O}\kern -.42em {\it `}\kern + .20em}}
\def\RR{{{\rm l}\kern - .15em {\rm R} }}
\def\PP{{{\rm l}\kern - .15em {\rm P} }}
\def\L2{{{\sf L}^2}}
\def\H1{{{\sf H}^1}}
\def\PN2{{\PP_{N}-\PP_{N-2}}}

\def\complex{{{\rm C} \kern - .53em {\rm l} \kern + .38em}}

\def\a1{{ | \lambda_{\min} |}}

\def\l1{{   \lambda_{\min}  }}

\def\bu0{{\underline {\bf 0}}}

\def\bu{{\bf u}}

\def\bx{{\bf x}}
\def\bux{{\underline {\bf x}}}

\def\Oh{{\hat \Omega}}

\def\ub{{\underline b}}

\def\ue{{\underline e}}

\def\up{{\underline p}}

\def\ur{{\underline r}}

\def\uu{{\underline u}}

\def\uw{{\underline w}}
\def\ux{{\underline x}}

\def\uz{{\underline z}}
\def\uzb{{\bar {\underline z}}}

\def\u0{{\underline 0}}
\def\1u{{\underline 1}}

\def\tR{{\! \widetilde R}}

\def\tA{{\tilde A}}

\begin{document}
\let\WriteBookmarks\relax
\def\floatpagepagefraction{1}
\def\textpagefraction{.001}

\shorttitle{Fast Coarse Solvers}

\shortauthors{T. Ratnayaka, P. Fischer, L. Olson}

\title[mode=title]{Coarse Solvers for Exascale Solution of Poisson Problems}



\author[1]{Thilina Ratnayaka}
\ead{rbr2@illinois.edu}
\author[2]{Paul Fischer}
\ead{fischerp@illinois.edu}
\author[3]{Luke Olson}
\ead{lukeo@illinois.edu}
%
\credit{software, analysis}

\cortext[1]{Thilina Ratnayaka}



\begin{abstract}
We present a two level Schwarz method as an alternative to Algebraic Multigrid
method (AMG) used as the last level (coarse) solver of the $p$-multigrid $p$MG
preconditioner for pressure Poission equation resulting from Spectral/Finite
element descretization of incompressible Navier-Stokes eqaution.
Proposed Schwarz method consits of a local problem in the original $p$MG coarse
space and a global coarse problem.
Main contribution of the paper is a novel, structured and a nonnested coarse
space for the global coarse problem.
Structured nature of the proposed global coarse space enable communication-free
interpolation between the original $p$-multgrid coarse space and the global coarse problem.
We demonstrate the effectiveness of the proposed method compared to the state of
the art AMG solver BoomerAMG by a series of experiments performed using Nek5000/RS,
a suite of highly scalable incompressible Navier-Stokes solvers, on Summit/Frontier
supercomputers at Oak Ridge Leadership Computing Facility.
\end{abstract}




\maketitle


 \section{Introduction}

$p$-multigrid ($p$MG) is an effective preconditioning strategy
for the scalable solution of elliptic problems that are discretized by high-order
finite-element (FEM) or spectral-element methods (SEM)~\cite{lottes05, mfem21,
hiptmair2002, malachi2022a}.
For example, using a local polynomial order $N=7$ with the SEM,
a multigrid smooth/restrict/prolongate schedule of polynomial orders
$N=7\to5\to3\to1$ is effective and reduces the number of degrees-of-freedom (dofs)
per element in 3D from $343\to125\to27\to1$.  Through local (fast) tensor-product
operations that require only $O(N^3)$ memory references and $O(N^4)$ work
per element, the method is highly efficient~\cite{dfm02}.   While this communication-minimal relaxation
process leads to a two-order-of-magnitude reduction in problem size, the
resulting coarse grid problem for $N=1$, denoted $A_c \uu_c = \ub_c$, is communication
intensive because $A_c^{-1}$ is \textit{completely dense}.  This is expected as the
Green's functions for the Poisson problem have slow decay; as a result,
data in any entry of the distributed coarse input, $\ub_c$, has an
immediate impact on \textit{all} values of the distributed solution, $\uu_c$,
which implies that the communication is all-to-all.

In this paper, we consider development of coarse grid solvers tailored to exascale
architectures (specifically, GPU-based platforms with high node counts).  The
target application is the Poisson problem that arises when simulating unsteady
incompressible flow with spatial discretization based on the SEM~\cite{dfm02,pat84}.
The SEM is a heterogeneous discretization that features
\textit{locally structured} degrees-of-freedom $u^e_{ijk}$,
$(i,j,k)\in [0,\dots N-1]^3$, which are ordered lexicographically within each of $E$
elements, $\Omega^e$, $1 \le e \le E$.  To accommodate complex geometries, the computational
mesh comprising the elements is \textit{globally unstructured}.  The mesh
restrictions are simply that each $\Omega^e$ should be an invertible (and
well-conditioned) map of the reference element $\Oh := [-1,1]^3$, that the
elements are nonoverlapping,\footnote{Overlapping SEM meshes in the spirit of
overset methods are, however, feasible~\cite{merrill16,mittal20b}.} and that
their union covers the computational domain; that is, $\Omega = \cup_e
\Omega^e$.
For $p$MG~\cite{malachi2022a} or multilevel Schwarz~\cite{fischer04,lottes05}
preconditioners, the heterogeneous nature of the SEM leads to a natural decomposition
of structured (tensor-product) work on the local mesh, which accounts for the majority
of the floating-point operations (flops), and communication-intensive operations on
the unstructured element-vertex (i.e., $N=1$) mesh, referred to as the \textit{coarse grid}.

It is well known that the \textit{coarse grid solve}\footnote{We
   refer to the coarse grid \textit{solve} throughout the text, but
   the coarse grid problem is rarely solved exactly because the level-$N$ outer iteration is
   embedded in a Krylov-subspace projection (KSP) such as conjugate gradients
   or GMRES\@.  For the coarse grid, a simple relaxation or other
   error-reduction scheme often suffices.}
is challenging in a parallel setting when the number of processes,
$P$, is large.
It has been observed
that the coarse grid solve is the only operation
in the solution of elliptic partial differential equations (PDEs) whose cost
\textit{increases} with $P$~\cite{fischer15}.  The other terms scale like
$(n/P)^\gamma$ for some $\gamma \geq 0$, where $n$ is the total number of dofs.
While terms with $\gamma < 1$ will impede ideal speed-up in the case of strong-scaling,
where $P$ is increased for fixed $n$, it is only the coarse grid solve that impacts weak
scaling, where $P$ is increased while the number of dofs local to each process, $n/P$,
is fixed.  As $P$ increases, the coarse grid solve
invariably sets the limits on strong-scaling, which directly impacts
time-to-solution in many applications.

The importance of the coarse grid solve for exascale applications is
illustrated in~\citet{sc22}, where it accounts for 45\% of the flow simulation
time when running on $P=\num{27648}$ Nvidia V100 GPUs on the Summit
supercomputer at the Oak Ridge Leadership Computing Facility (OLCF).
In this case, the Navier-Stokes problem has $E=98$ million elements of
order $N=8$ (leading to $n=51$ billion dofs) and the coarse solve in the
$p$MG-preconditioned pressure Poisson problem uses a single $V$-cycle of
BoomerAMG~\cite{boomer-amg-2002}.  For $p$MG on tensor-product
elements, the coarse grid size is $n_c \approx E$.  Current exa- and
pre-exascale platforms are routinely enabling $E=10^8$--$10^9$, resulting in
large ``coarse'' systems.  In~\citep{sc22}, the coarse grid overhead is
indrectly reduced by adding more smoothing sweeps at the finer $p$MG levels,
thereby reducing the number of outer GMRES iterations and
consequent visits to the $N=1$ level of $p$MG.

In this paper, we consider a direct approach to reducing communication for the
$N=1$ coarse solve $A_c \uu_c=\ub_c$ by introducing a low-communication
two-level overlapping Schwarz
smoother with a novel nonnested coarse space, as an alternative to AMG\@.  The
overlapping systems (size $\approx E/P$) are highly localized and well suited to computation
on the GPU, yielding full $P$-fold parallelism with minimal off-device data
transfer.  The nonnested (reduced) coarse space requires only a few (e.g.,
4--10) dofs per process and the corresponding reduced system, $A_r \uu_r = \ub_r$
(size $n_r =O(P) \ll E$), is solved using the communication-minimal $XX^T$
algorithm developed in~\cite{fisc96,tufofisc01}.

Rest of the paper is organized as follows.
In Section~\ref{sec:rev}, we review past coarse grid solve strategies
and underscore the need for a novel approach in the context of GPU-based
exascale platforms.
In Section~\ref{sec:schwarz}, we present the proposed two-level Schwarz
method, including a coarser system, $A_b \uu_b = \ub_b$, which is based
on a simple nonnested approximation space.
We discuss the background and details of the nonnested coarse space
in Section~\ref{sec:nnest}.
Numerical results with both Nek5000 and NekRS using four different meshes
are shown in Section~\ref{sec:results}.
Conclusions and future directions are given in Section~\ref{sec:future}.

 \section{Survey of Coarse Grid Solvers}\label{sec:rev}



Solving the coarse problem using the explicit inverse of $A_c$ becomes
prohibitive in 3D for large $n_c$ due to the storage cost of $A_c^{-1}$.
Even if $A_c^{-1}$ is distributed across the number of processes $P$, this
method has a storage cost of $n_c^2/P$.
However, the fact that we can simply solve the system using $u = A_c^{-1}b$
is quite appealing since all the dot products can be performed in parallel.
Approach described in~\cite{alvarado1992highly} try to exploit this feature by
inverting the Cholesky factor $L$ of $A_c$ using a partitioned inverse approach.
Any unit lower triangular matrix of size $n_c$ can be expressed as a product of
$n_c$ elementary matrices $L_i$ i.e., $L = \prod_{i=1}^{n_c} L_i$, $L_i$ is of
the form $I + m_ie_i^T$ where $e_i$ is the $i^{th}$ unit vector and $m_i$ is a
vector with first $i$ values being zeros.
Then these $L_i$ elementary matrices are grouped together to form $k$ factors
such that $L = \prod_{i=1}^k P_i$ where each factor $P_i$ has the property that
$P_i^{-1}$ can be represented using the same space as $P_i$.
The solution $\uu = L^{-1}\ub = \prod_{i=k}^1P_i^{-1}\ub$ can be calculated by
doing $k$ matrix vector multiplications.
Communication for this approach scales as $O((\log_2{P})^2)$ and is
not optimal.

During time stepping of a computational fluid dynamics simulation, system
$A_c\uu_c=\ub_c$ is solved multiple times with different right hand sides for
the same constant matrix $A_c$.  In time transient problems, successive right
hand sides (RHS) often share enough information so a good initial guess based
on previously generated Krylov subspaces can yield significant reductions in
the global solve costs \cite{fisc93,fisc98}.  Farhat and Chen \cite{farhat94}
suggest using this approach to solve $A_c\uu_c = \ub_c$ for successive RHS by
using a slightly modified version of Conjugate Gradient (CG) iteration to make sure that
the search directions generated for new RHS are $A-$orthogonal to previously
generated Krylov subspace.  When a solution $\uu_c$ for a new RHS $\ub_c$ has
to be found, the initial guess is computed as the solution for a restricted
versions of $A_c \uu_c = \ub_c$ in the previously constructed Krylov subspace.
The dominant cost for this method is this projection step to calculate the
initial guess since the cost for subsequent CG iterations are negligible.
Each of the CG iterations following this initial projection incurs a latency
cost of $4\alpha\log_2 P$ due to the reductions associated with CG.

An alternative approach is to use the sparse direct factorization $A_c^{-1}=XX^T$
proposed in \cite{tufofisc01}.
Main idea is to find a set of $k$, $A_c$-conjugate sparse basis vectors
$X_k = [\ux_1\;\ux_2\;\dots \ux_k]$ satisfying $\ux_i^T A_c\ux_j=\delta_{ij}$,
where $\delta_{ij}$ is the Kronecker delta function.
Then an approximation to $\ux_c$, ${\bar \ux}_k$ can be calculated by projecting
$\ub_c$ in to this space i.e., ${\bar \ux}_k=X_kX_k^T\ub_c$.
If $k=n_c$, this approach yields the exact solution, $\ux_c={\bar \ux}_{n_c}$.
By using a nested-disection ordering, the authors showed that the communication
complexity for a system arising from a discrete Poisson problem in $d$ space
dimensions is $2(\alpha^*+\beta^* Cn_c^{(d-1)/d})\log_2 P$, where $C$
is a constant and $(\alpha^*,\beta^*)$ are interprocess latency
and inverse bandwidth respectively.
For sufficiently small systems this cost is $\approx 2\alpha^* \log_2 P$, which
is near optimal.


 \section{Two-Level Coarse Solver}\label{sec:schwarz}

\noindent
In this section, we develop a low-communication preconditioner, $M_c$, to
approximate the coarse-level ($N=1$)  system matrix, $A_c$.
Given that the full pMG system is embedded in a Krylov subspace projection
(KSP), an exact solution for the $N$=1 system is not required; a single sweep
of a multilevel preconditioner is often sufficient.  A potential consequence of
an inexact solve, however, is an increase in the number of KSP iterations, so
any economization must be weighed against the overall cost of the outer
KSP solver, as discussed in the examples of Section \ref{sec:results}.

We base $M_c$ on the two-level additive Schwarz method (ASM) of Dryja and
Widlund \cite{drywid87,drywid89,drywid94}.  To simplify notation, we drop the
subscript $c$ for the coarse problem and simply refer to the system as $ A\uu =
\ub$ and to the preconditioner as $M$, unless otherwise needed. In this
context, $A$ is assumed to be a sparse $n \times n$ symmetric positive definite
(SPD) matrix derived from an unstructured {\em nodal} finite element (FEM)
discretization of the Poisson equation with homogenous Dirichlet boundary
conditions.  (In our case, the $A_c$ system is based on a Galerkin restriction
of the originating order-$N$ SEM operator.) Many features of this two-level
method would work equally well for other PDEs, however, and for other
nodal-based discretizations such as finite differences or finite volumes.  We
reiterate that our target problem size will involve only $n/P \approx
10^3$--$10^4$ dofs per GPU and that this problem is thus intrinsically
communication bound, despite the potentiality that $n$ might be as large as
$10^9$.

To set the stage for later developments, we introduce two representations for
FEM basis functions on the domain, $\Omega$, having boundary $\dO = \dO_N
\cup \dO_D$.  We assume that the boundary conditions for the PDE are
homogeneous Neumann on $\dO_N$ and homogeneous Dirichlet on $\dO_D$.
Let $X^h_0 =$span$\{\phi_j(\bx)\}$ be the set of $C^0$ continuous basis
functions satisfying $\phi_j(\bx_i)=\delta_{ij}$ (the Kronecker delta function)
that forms the FEM basis on nodal points $\bx_i$.  We assume that all basis
functions vanish on $\dO_D$ and that any trial solution $u(\bx) \in X^h_0$ can
be written as
\begin{eqnarray} \label{eq:glob}
   u(\bx) &=& \sum_{j=1}^n \,  \phi_j(\bx) \, u_j
   \;=\; \sum_{e=1}^E \sum_{k=1}^{n_v} \,  l^e_k(\bx) \, u^e_k.
\end{eqnarray}
The first expression on the right is the {\em nodal} (or {\em global}) form
involving unknown basis coefficients, $u_j$, often referred to as global
degrees-of-freedom (dofs).  The second expression is the {\em elemental} (or
{\em local}) form.  We assume that $\Omega=\cup_e \Omega^e$ represents a
decomposition of $\Omega$ into $E$ nonoverlapping trianglular or quadrilateral
elements (tetrahedral or hexahedral in 3D) of size $O(h)$, each having $n_v$
vertices, $\bx_k^e$.  For each element, there are $n_v$ local Lagrange interpolants,
$l_k^e(\bx)$, $k=1,\dots,n_v$, which vanish outside of $\Omega^e$.
There are a total of $n_l = E \cdot n_v$ local
basis coefficients, $u^e_k$.  To impose $C^0$ continuity on $u(\bx)$, we
additionally have a constraint on these coefficients, $\bx_k^e = \bx_{k'}^{e'}
\, \Longrightarrow \, u_k^e = u_{k'}^{e'}.$ The standard way to impose this
constraint is to start with the global representation and to map the
coefficients $u_j$ to their local counterparts, $\bx_k^e = \bx_j \,
\Longrightarrow \, u_k^e = u_j.$ The transformation from global to local
coefficients can be expressed as a matrix-vector product, $\uu_L=Q \uu$, where
$Q$ is an $n_l \times n$ Boolean matrix, $\uu_L = \{ u_k^e \}$ is the set of
local unknowns, and $\uu=\{ u_j \}$ is the set of global unknowns.  Application
of $Q$ is a global-to-local copy operation, whereas application of $Q^T$
involves summation (contraction) of local basis coefficients to their
global counterparts.

The FEM discretization of the Poisson problem, $-\nabla^2 u = f$, is based on the
weak form.  Let $(v,u):=\int_{\Omega} v\,u\,dV$ be the $L^2$ inner-product on
$\Omega$, perhaps approximated by a suitable quadrature rule.  Then the
standard Galerkin projection approach leads to the FEM system matrix, $A$, having
entries $a_{ij} = (\nabla \phi_i,\nabla \phi_j)$, and right-hand side, $\ub$,
with entries $b_i=(\phi_i,f)$. It is easy to show that $A=Q^T A_L Q$, where
$A_L$=block-diag($A^e$) comprises local system matrices having entries
\begin{eqnarray}
  a_{ij}^e = \int_{\Omega^e} \nabla l_i^e \cdot \nabla l_j^e \, dV
\end{eqnarray}
Application of $Q^T$ and $Q$ to $A_L$ is referred to as the matrix assembly
process, which we will use later in the development of our reduced coarse-space
operator.

For applications on $P$ distributed-memory compute units (i.e., processes or
MPI ranks), it is natural to cluster the elements into $P$ contiguous
subdomains using, say, recursive spectral bisection \cite{pothen90} or other
partitioning strategy that minimizes the number of shared dofs on the subdomain
interfaces.  We denote these subdomains as $\Omega_p$, $p=1,\dots,P$, and map
data and solution values associated with $\Omega_p$ to process $p$.
Nodal data on $\dO_p$, the boundary of $\Omega_p$, is redundantly represented
on any process that shares those boundary nodes.  With this element-based
partition, nearest-neighbor communication in our implementation typically
arises when computing matrix-vector products of the form, $\uw=A \uu$.
In local form, this product is expressed as $\ur_L = A_L \uu_L$, which is
communication free because $A_L$ is block-diagonal, followed by
direct-stiffness summation, $\uw_L = QQ^T \ur_L$, in which shared nodal values
are summed by $Q^T$ and redistributed by $Q$ \cite{dfm02}.

For the Schwarz method, we extend each subdomain $\Omega_p$ by one or more
layers of elements from neighboring processes to form a set of {\em
overlapping} subdomains, $\Ol_p$, with overlap $\delta = O(h)$.  For each
domain, we define a Boolean restriction operator, $R_p$, such that $\uu_p :=
R_p \uu$ returns the vector of nodal values that are {\em interior} to $\Ol_p$
(i.e., excluding values on $\partial \Ol_p$).  Note that $\uw_p := R^T_p \uu_p$
extends (by zero) a vector of local subdomain values on $\Ol_p$ to a global
vector of length $n$.  With this decomposition, the additive Schwarz
preconditioner is
\begin{eqnarray} \label{eq:asm}
\uz  &=&
\sum_{p=1}^P \, \tR_p^T A_p^{-1} R_p \ur
\;+\; J A_r^{-1} J^T \ur ,
\end{eqnarray}
which is the sum of local subdomain operators (subscript $p$) and a reduced
coarse grid operator (subscript $r$), which we introduce later.  The method is
naturally parallel because the subdomain problems can be solved independently.
The reduced coarse-space system, $A_r := J^T A J$, with $J$ a coarse-to-fine
interpolation matrix, is not trivially solved in parallel but it is of modest
size compared to $A$ and can be solved using an $XX^T$ factorization, as we
discuss below.

The local systems are $A_p := R_p A R_p^T$, which constitute the principal 
submatrices of $A$ that correspond to the interior nodes of $\Ol_p$.  In
(\ref{eq:asm}), we also introduce the prolongation operator, $\tR_p^T$, which
is equal to $R^T_p$ for the classic ASM.  Once the vectors are formally
extended by zero, they can be added.  Thus, summation of the local solutions,
$\uzb_p:=R_p^T A_p^{-1} R_p \ur$, amounts to summing shared nodal values in
regions of overlap.   These local solution values must be communicated so that
each neighboring process can complete the sum for all nodes in $\Omega_p$.
Alternatively, one can define $\tR^T_p$ such that each process ignores all
nonlocal components from the local Schwarz solve. This approach yields the {\em
restricted additive Schwarz} (RAS) method of Cai and Sarkis \cite{ras99}, which
is in some cases more rapidly convergent than the traditional ASM.  For RAS,
$\tR_p^T$ is the same as $R_p^T$, save that columns associated with $\Ol_p
\backslash \Omega_p$ are null, so less communication is required.  One still
needs communication to enforce continuity for the redundantly-stored surface
nodes, however.

It is clear that the local solutions of the Schwarz method require a minimal
amount of communication.  Only surface values need to be communicated
(particularly in the minimal-overlap case).  On GPUs, one can solve (exactly or
approximately) the local problem on the device.  There is no need to
communicate $O(n/P)$ data between the device and the host.   Unlike multigrid,
two-level Schwarz incurs nearest neighbor communication only on the fine level.
The remaining communication is deferred to the coarse problem.


As is well known, Schwarz methods must be augmented with a coarse space
approximation to realize convergence rates that are independent of $P$ (i.e.,
that are scalable) \cite{sbg96}.  The most common approaches are the additive
variant (\ref{eq:asm}) or a hybrid method that has a locally additive step,
\begin{eqnarray} \label{eq:asm1}
\mbox{ASM$_1$:} \;\;\;\; \uz_{\mbox{\footnotesize \em loc}}
&=& \sum_{p=1}^P \, \tR_p^T A_p^{-1} R_p \, \ur,
\end{eqnarray}
followed by a coarse-space correction,
\begin{eqnarray} \label{eq:asm2}
\mbox{ASM$_2$:} \;\;\;\; \uz \;\; &=&
\uz_{\mbox{\footnotesize \em loc}}
 \;+\; J A_r^{-1} J^T \, (\ur\,-\, A \uz_{\mbox{\footnotesize \em loc}}).
\end{eqnarray}
The multiplicative approach (\ref{eq:asm1})--(\ref{eq:asm2}) serializes the
local and coarse solve steps but is often more effective than (\ref{eq:asm}) in
reducing the overall number of iterations in the full system.

In (\ref{eq:asm}) and (\ref{eq:asm2}), $J$ is an interpolation matrix from the
reduced coarse-space basis functions to the nodal points, $\bx_j$.  The main idea
is that the interpolated range of $J$ should provide an efficient representation
of the low-wave number components of the Poisson operator on $\Omega$.

\noindent
For nested coarse spaces with exact local and coarse solves, the condition
number for the preconditioned system is
\begin{eqnarray}
\kappa(M^{-1}\!A) &=& O(1+H/\delta),
\end{eqnarray}
where $\delta$ is the amount of domain overlap and $H$ is the diameter of the
subdomains, which is assumed to be the same as the diameter of the coarse-space
elements.  Cai \cite{cai94} gives an alternative estimate for the nonnested
case,\\[1ex]
\noindent

\noindent
where $H$ a characteristic size of the coarse-space elements (i.e., boxes, in
the notation of the next section), which is independent of the subdomain size.
Under slightly different constraints on the coarse space, Chan, Smith, and Zou
\cite{chan1996overlapping} give the following bound for additive Schwarz with
nonnested spaces,
\begin{eqnarray} \label{eq:csz_eq21} 
\kappa(M^{-1}A) & \leq & C \left(1+\frac{H}{\delta} \right)^2.
\end{eqnarray}
The authors remark (Remark 5 \cite{chan1996overlapping}) that this bound can
be improved to
\begin{eqnarray} \label{eq:csz_rem5} 
\kappa(M^{-1}A) & \leq & C \left(1+\frac{H_{\max}}{\delta} \right),
\end{eqnarray}
if the subdomains $\Omega_p$ form a quasi-uniform triangulation of $\Omega$
and if $H \leq \beta H_{\max}$ for some fixed constant $\beta$ where
$H_{\max} = \max \mbox{diam} \Omega_p$.


\section{Nonnested Coarse Spaces}\label{sec:nnest}

As noted earlier, a coarse-space is essential for efficient iterative solution
of large-scale problems.  The role of the reduced coarse space is to provide a
mechanism for generating long wavelength approximations to the solution with
relatively few degrees of freedom.  Without such an approximation, many
iterations may be required for low-wavenumber components of the solution to
emerge in the iteration process, particularly if one is using restarted-GMRES
rather than a projection-based iteration.  For unstructured mesh problems, it
is not immediately obvious how to construct a coarse space that utilizes the
underlying approximation space, $X^h_0$.  Aggregation (either smoothed or
unsmoothed)  and fine-coarse (FC) splittings are two of the most popular
approaches in AMG.  Nonnested coarse spaces have been discussed in both
multigrid and Schwarz contexts~\cite{cai94,bramble1991,xu1996auxiliary}.  Here, we follow the nonnested
approach as it affords significant flexibility and the desired level of
granularity for our targeted exascale platforms.  Moreover, we can leverage the
communication-minimal $XX^T$ coarse grid solver that is designed for problems
at this scale \cite{fisc96,tufofisc01}.

Without the restriction of nested spaces, one can develop a reduced coarse
space $X^r$ based on relatively simple interpolants.  Let
$\hX^r=$span$\{\Phi_j(\bx)\}$, $j=1,\dots,n_r,$ define a set reduced-space
interpolation functions on $\Omega^r \supset \Omega$ and define the
interpolation matrix $J$ as having entries $J_{ij} := \Phi_j(\bx_i)$.   Then
the reduced approximation space is given by $X^r=$span$\{\Psi_i(\bx)\} \subset
X_0^h$, where
\begin{eqnarray}  \label{eq:Psi}
  \Psi_j(\bx) &=& \sum_{i=1}^n \, \phi_i(\bx) J_{ij}. 
\end{eqnarray}
Note that for any nodal point $\bx_i$ in the originating FEM mesh we have
$\Psi_j(\bx_i) = J_{ij} = \Phi_j(\bx_i)$.  However, the functions $\Psi_j$ are
properly in $X_0^h$---they satisfy the original continuity requriments and
boundary conditions---whereas the $\Phi_j$s are not so constrained.
We can choose the basis for $\hX^r$ based solely on efficiency and
ease-of-implementation considerations.  In fact, for advection-diffusion
problems there are stability advantages to admitting discontinuous bases in
$\hX^r$, as we illustrate in the Appendix.  Discontinuous bases for $\hX^r$ are
also attractive for local mesh refinement in the reduced space as one requires
no special treatment for hanging nodes.

With the definition (\ref{eq:Psi}), the Galerkin statement for the Poisson
problem in the reduced space is, {\em Find $u_r(\bx) \in X^r \subset X_0^h$
such that}
\begin{eqnarray}
(\nabla v,\nabla u_r) &=& (\nabla v, \nabla u) \;\; \forall \; v \, \in X^r.
\end{eqnarray}
In matrix form, the reduced-space basis coefficients
$\uu_r=[u_{r,1}\;\dots u_{r,n^r}]^T$ are governed by the reduced system,
$A_r \uu_r = J^T \ub$, with the SPD system matrix, $A_r:=J^T A J$.

\begin{figure}[t] \centering
{\setlength{\unitlength}{1.0in}
   \begin{picture}(6.500,1.90)(0.03,.8)
     \put(0,0){\begin{picture}(0,0)(0,0)
       \put(0.00,0.00){\includegraphics[width=2.25in]{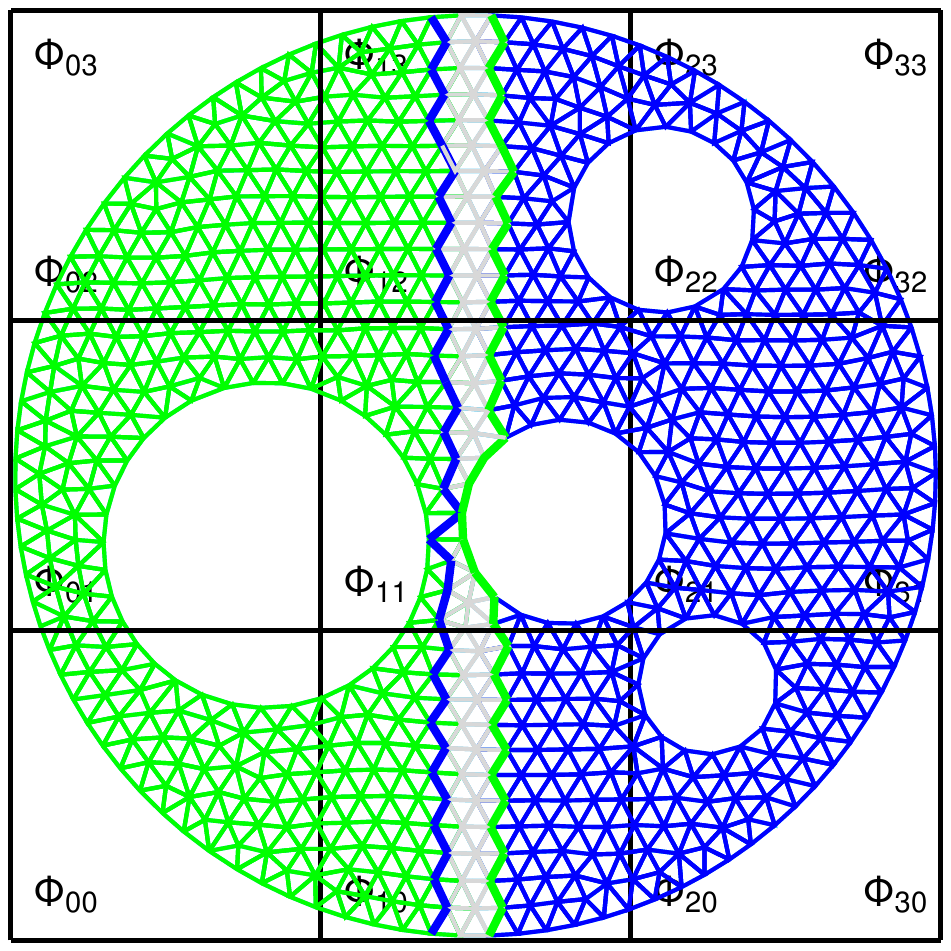}}
       \put(2.20,0.00){\includegraphics[width=2.25in]{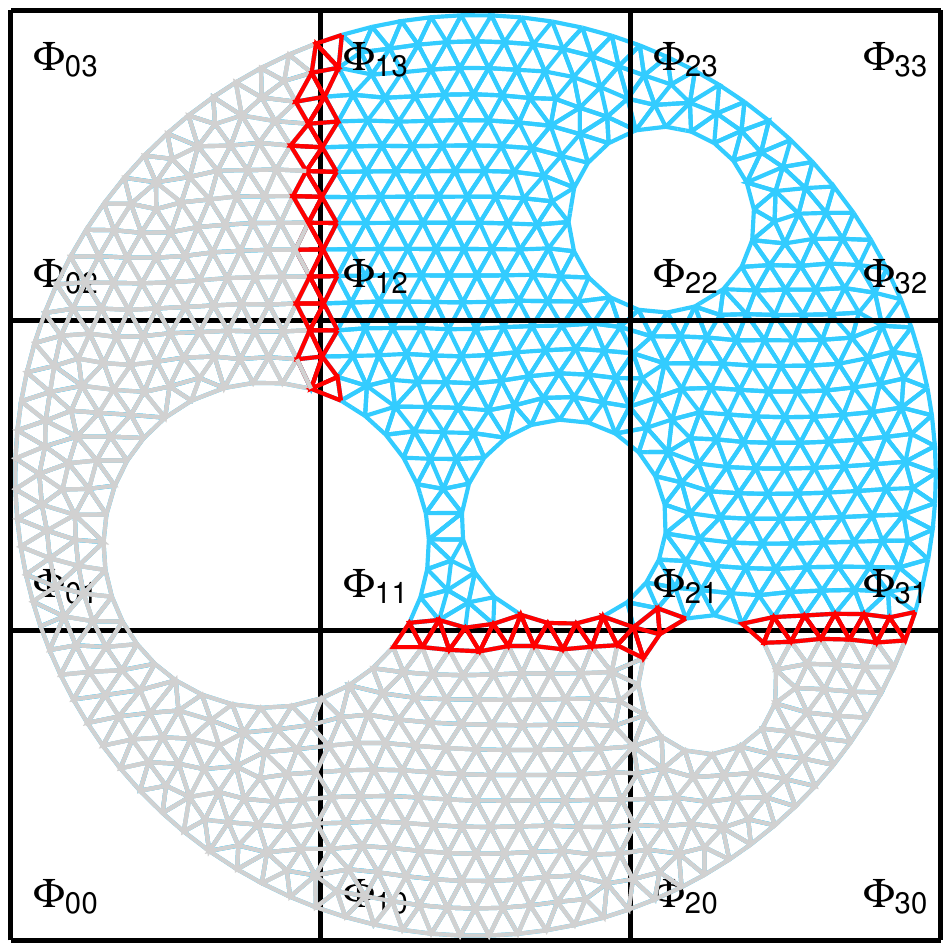}}
       \put(4.40,0.00){\includegraphics[width=2.25in]{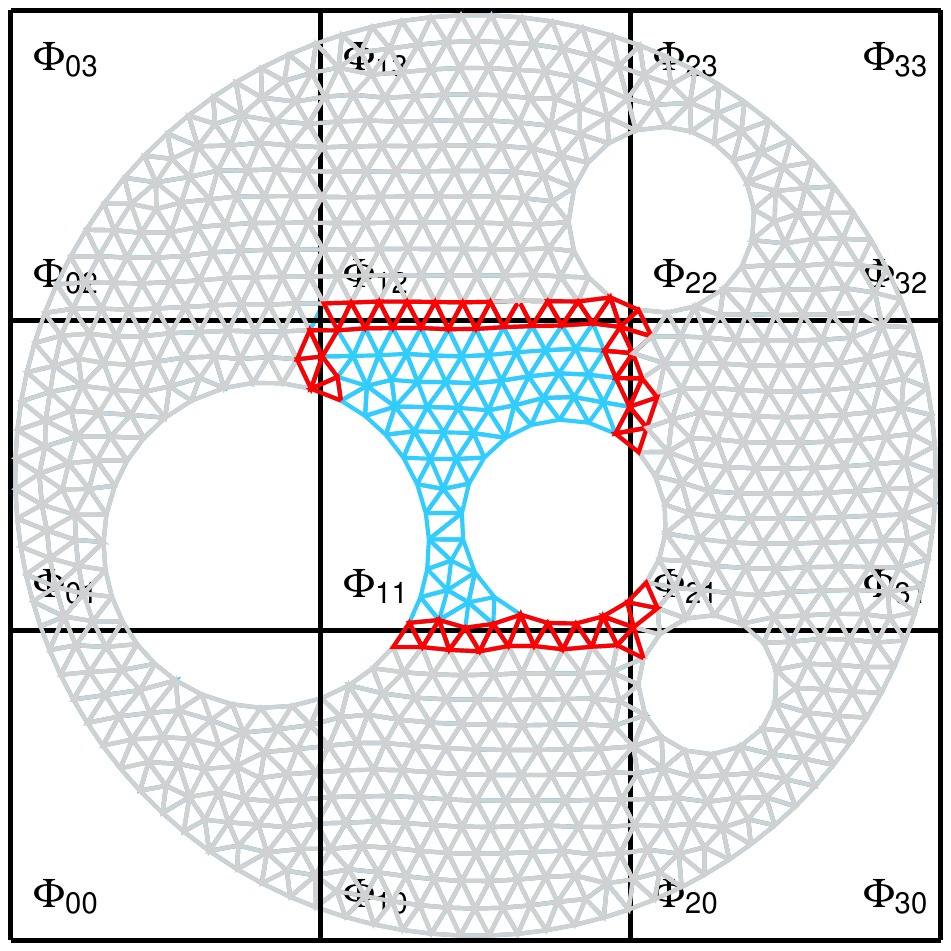}}
     \end{picture}}
     \put(0.1,2.2){\begin{picture}(0,0)(0,0)
       \put(0.00,0.00){(a)}
       \put(2.20,0.00){(b)}
       \put(4.40,0.00){(c)}
     \end{picture}}
   \end{picture}
}
\caption{\small
Two-level Schwarz illustration.
%
(a) Parallel partition ($P=2$) of the $N=1$ mesh: each rank solves a local
Poisson problem on their respective shaded region (green or blue), including a
one- or two-element overlap extension, shown in gray.
(b) Vertex-based support of the coarse grid interpolants: shown in blue
is the support for $\Phi_{22}$.  The red triangles indicate border elements, for
which support of $\Psi_{22}$ would interact with $\Psi_{0k}$ and $\Psi_{k0}$,
$k=1$ and 2, which would in general give rise to a $5 \times 5$ stencil for
these lexicographically-ordered coarse-space functions.
(c) Element-based illustration of the coarse decomposition:  if the
(red) overlapping elements are eliminated then the coarse stencil will
only be $3 \times 3$ and the coarse-space operator, $A_r$, can be formed
using standard FEM assembly techniques.
\label{fig:support}
}
\end{figure}

Here, we take $\hX^r$ to be the space of tensor products of piecewise linear
interpolants that cover a rectangular box in $\RR^d$.  The 1D linear bases give
rise to box-like elements in 2D and 3D, as illustrated in Fig.
\ref{fig:support}.  With this basis, it is easy to construct a structured
coarse-space to interpolate to the fine nodes.  Information regarding topology
or boundary conditions associated with $\Omega$ is not required.  Aside from
some solvability constraints (discussed below), the domain of $X^r$ can be a
simple superset of $\Omega$, which implies that minimal geometric information
is required to construct a viable set of interpolants, $\Phi_j(\bx)$.  We
simply find the limits of $\Omega$ in each of the principal directions and
partition the space into a set of $n_b = b_x \times b_y (\times b_z)$ boxes
that span these dimensions.  We choose the number of boxes to be a small
multiple of $P$ (e.g., 4--10) so that there are just a few dofs per process
and to ensure that the box size, $H$, is greater than the characteristic finite
element size, $h$.  The number of reduced-space unknowns is $n_r \leq
(b_x+1)(b_y+1)(b_z+1)$, with the upper-bound corresponding to the total number
of box vertices.  Inequality holds for cases when there are empty boxes or,
more precisely, when a larger value for $n_r$ would yield linearly dependent
columns in $J$.

In the example of Fig. \ref{fig:support}, the reduced-space problem has 16
unknowns, one for each vertex associated with the coarse bilinear interpolants,
here denoted with a dual subscript, $\Phi_{ij}$, $i,j \in \{0,\dots,3\}^2$.
  Figure \ref{fig:support}(b) illustates the vertex-oriented interpretation of
the basis functions, with the support for $\Phi_{22}$ shown in blue.  The
support for $\Psi_{22}$, which includes the blue and the red triangles, can be
seen to extend beyond the squares that surround $\Phi_{22}$.
  Figure \ref{fig:support}(c) gives an ``element'' based interpretation of the
support for the center-most box.  Here, every triangular element marked in blue
is influenced only by the four basis functions, $\Phi_{11}$, $\Phi_{21}$,
$\Phi_{12}$, and $\Phi_{22}$.  The red triangles indicate border elements whose
support spans both the blue region and the adjacent boxes.  A consequence of
these overlapping elements is to increase the overall support of the individual
basis functions, $\Psi_{ij}(\bx)$, which leads to reduced sparsity in $A_r$.
We address this issue in the next subsection.  Note that there is no need to
integrate or differentiate the interpolants $\Phi_{ij}$. All of the physics is
contained in $A$ and consequently embedded in $A_r := J^T A J$.  If $J$ is of
full rank and $A$ is SPD, then $A_r$ will be SPD and the {\em reduced coarse
problem}, $\uz_r = J A_r^{-1} J^T \ur,$ will be solvable.  This system will
provide the best $A$-norm approximation to the solution $A^{-1} \ur$ in the
range of $J$.  By construction, $\uz_r$ will provide the desired
long-wavelength approximation to the unknown error.

Note that the holes in the domain of Fig. \ref{fig:support}
pose no particular difficulty.  If the
boundaries of the holes are part of $\dO_D$ the solution will vanish on
those boundaries and the reduced approximation, $\uz_r$, will respond
accordingly.  Similarly, if the boundary for the hole is part of $\dO_N$ the
interpolated solution will float according to the nontrivial values in $\uz_r$
associated with nodes $\bx_i \in \dO_N$.
  A significant attraction of the nonnested spaces is that one can generate a
simple coarse space that can be applied without detailed attention to geometry
or boundary conditions.  Situations where concern for geometry is warranted,
however, include external flows such as flow over aircraft.  In these cases,
one might find that the far-field elements are larger than $H$, which can
result in boxes that are interior to $\Omega$ but which contain no fine-scale
vertices, potentially leading to null columns in $J$ since those bases
would have no target interpolation points.  Unlike vertices exterior
to $\Omega$, interior $\Phi_j$s cannot be trivially discarded because they are
needed to formally provide function continuity across $\Omega$.  In such
situations it is necessary to coarsen the reduced space and perhaps augment it
with local mesh refinement in the near-field regions.  The flexibility of
the current approach provides a significant benefit in this case because
one does not need continuity of $\Phi_j$ at hanging nodes.

In a high-level language like Matlab, implementation of the reduced space
solver is straightforward.  One begins with a sparse matrix $n \times n_r$
matrix $J$ that is void and then marches through each nodal point, $\bx_k$,
$k=1,\dots,n$, to identify its bounding box.
  Given the bounding box, one can identify each of the $2^d$ box vertices and
construct the linear interpolants that map the box vertex values to the node in
question.
  Consider for example a 2D case with a $b_x \times b_y$ array of boxes.  Let the
columns of $J$, which represent the discrete support of $\Phi_{ij}$, be ordered
lexicographically.  Then
\begin{eqnarray}
J &=& [ \Phi_{00}(\bux) \; \;\; \Phi_{10}(\bux) \; \;\; \cdots \; \;\;\Phi_{b_x b_y}(\bux) ],
\end{eqnarray}
where $\bux = [ \bx_1 \;\; \bx_2 \;\; \cdots \;\; \bx_n]^T$ is the vector of
nodal points.  Because the support of each $\Phi_{ij}$ is compact, $J$ will be
sparse with at most $2^d n$ nonzeros for a $d$-dimensional problem.  (Each
nodal point $\bx_i$ is influenced by the corner vertices of its bounding box,
so each row of $J$ will have  $2^d$ nonzeros.) Construction of $J$ proceeds by
finding and evaluating the four bounding interpolants for each $\bx_i$.
One can flag the vertices in each box that has a nontrivial entry.  Any
bounding vertex that is not flagged is compressed out of $J$.
In our production $XX^T$ solver that is used to solve the $A_r$ system, the
setup phase does not require contiguous index sets, so columns that are not
passed into the setup call are automatically suppressed.

\begin{figure}[t]
\centering
{\setlength{\unitlength}{1.0in} \begin{picture}(6.600,1.900)(0.0,-.1)

 \put(0.0,0.2){\begin{picture}(2.2,1.3)(0,0)
   \put(0.0,0){\includegraphics[width=2.00in]
                                   {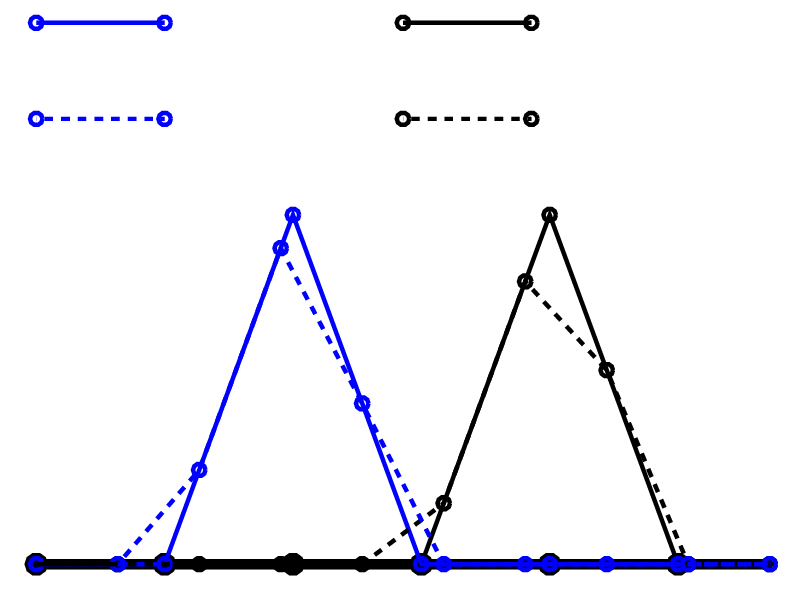}}
   \put(0.1,0.2){\normalsize (a)}
   \put(0.29,0.05)
    {$\underbrace{\hspace*{.83in}}_{\mbox{\normalsize $S(\Psi_{j-1})$}}$}
   \put(0.91,0.05)
    {$\underbrace{\hspace*{.83in}}_{\mbox{\normalsize $S(\Psi_{j+1})$}}$}

   \put(0.49,1.38){$\Phi_{j-1}$}
   \put(1.41,1.38){$\Phi_{j+1}$}

   \put(0.49,1.15){$\Psi_{j-1}$}
   \put(1.41,1.15){$\Psi_{j+1}$}

 \end{picture}}

 \put(2.2,0.2){\begin{picture}(2.2,1.3)(0,0)
   \put(0.0,-.02){\includegraphics[height=1.49in,width=2.00in]{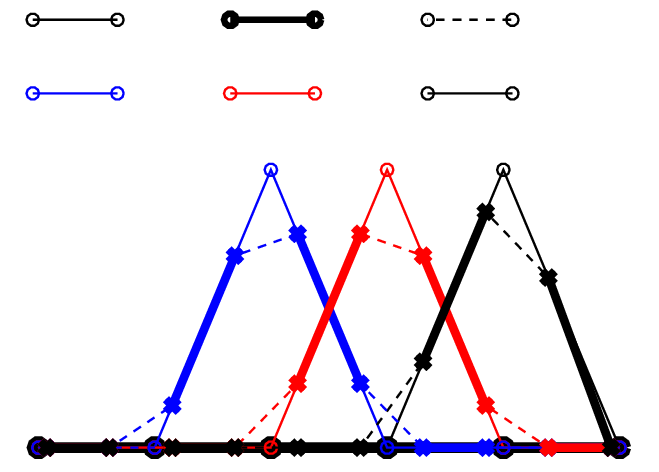}}
   \put(0.1,0.2){\normalsize (b)}

   \put(0.43,1.38){$\Phi$}
   \put(1.04,1.38){$\Psi$}
   \put(1.59,1.38){$\Psi_{\mbox{\tiny \em orig}}$}

   \put(0.43,1.15){${j}$-1}
   \put(1.04,1.15){${j}$}
   \put(1.59,1.15){${j}$+1}

 \end{picture}}

 \put(4.4,0.2){\begin{picture}(2.2,1.3)(0,0)
  \put(0,-.078){\includegraphics[height=1.55in,width=2.0in]{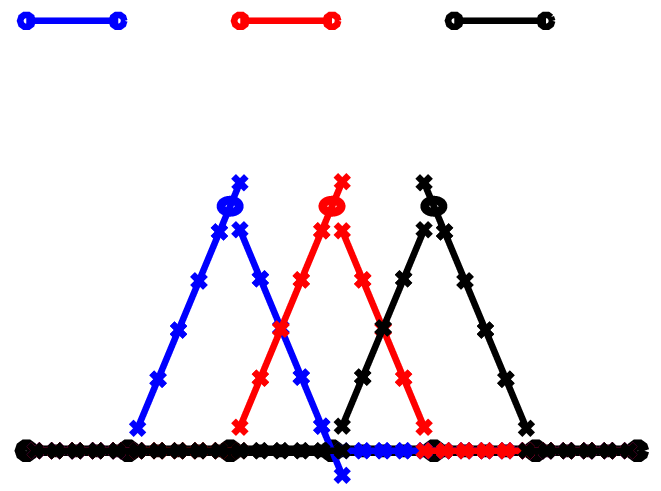}}
   \put(0.1,0.2){\normalsize (c)}
   \put(0.44,1.40){$\Psi_{j-1}$}
   \put(1.09,1.40){$\Psi_{j}$}
   \put(1.71,1.40){$\Psi_{j+1}$}
 \end{picture}}

\end{picture}}
\caption{
Reduced-space basis functions in 1D:
(a) standard, following Eq. (\ref{eq:Psi}), showing overlapping support
    of $\Psi_{j-1}$ and $\Psi_{j+1}$;
(b) gap-based support;
(c) element-centroid-based support.
The overlap in case (a) leads to a 5-point stencil while (b) and (c) 
yield a 3-point stencil such that $A_r$ is tridiagonal.
\label{fig:1D}
}
\end{figure}

It is helpful to recognize that the action of $J$ is to copy the reduced vertex
values to the $2^d$ cells in the box domain that surrounds each vertex.  That
vertex value is interpolated onto each FEM node within a given cell.  The
action of $J^T$, therefore, is to multiply the FEM nodal values by the
same weights and to then sum these onto the $2^d$ indvidual vertex values.
One can thus cast the formation of $A_r$ as a matrix
assembly process.  Consider the 2D case as an example and let $S_{ij}$
denote the cell subdomain covered by the support of
$\{ \Phi_{i-1,j-1} \Phi_{i,j-1} \Phi_{i-1,j} \Phi_{i,j} \}$.
Let $A_{ij}$ denote the $4 \times 4$ matrix associated with cell $S_{ij}$
that accounts for finite elements interior to $S_{ij}$.
(This matrix is the cell equivalent of $A^e$ introduced in the FEM context.)
For elements $\Omega^e$ that are interior to $S_{ij}$, one can form a local
$n_v \times 4$ interpolant, $J_{ij}^e$, that maps the four vertex values
for cell $S_{ij}$ to the $n_v$ nodes in $\Omega^e$.
Let $E_{ij}$ be the set of elements that are interior to $S_{ij}$ and
define
\begin{eqnarray} \label{eq:Aij}
  A_{ij} &=& \sum_{e \in E_{ij}} \, (J_{ij}^e)^T A^e J_{ij}^e.
\end{eqnarray}
In an ideal situation, one could simply assemble each of these $4 \times 4$
matrices to form $A_r$.  Unfortunately, this is not quite the case because the
procedure above applies only to elements {\em interior} to $S_{ij}$ (i.e., the
blue triangles in Fig. \ref{fig:support}(b)).  If we account for the red
triangles then the simple assembly procedure breaks down because of the
extended support of the $\Psi_{ij}$s.  Possible mitigation options are to
ignore the red triangles or to assign them decisvely to one cell or another
(e.g., by the location of their centroid).  These ideas are illustrated
in Fig. \ref{fig:1D}.

    Ignoring the red triangles is roughly equivalent to using an incomplete
quadrature rule that skips elements cut by cell boundaries.   The problem with
this approach is that one might have a chain of isolated triangles with
vertices on $\dO_N$, each cut by a cell boundary.  Representation of those
vertex values would be lost by the discard rule and the reduced solution would
not be propagated to those vertices.
  Assigning elements (as opposed to vertices) to cells, on the other hand has
several advantages.  First, all nodal values are represented.  Second, any cell
that has at least one element will be assigned at least $n_v$ rows in the
corresponding columns of $J$, thereby reducing the likelihood that columns of
$J$ will be linearly dependent.   One downside of the centroid assignment
approach is that the output is no longer continuous.  This situation can be
remedied, however, by applying $W Q_c Q_c^T$ to $\uz_r = J A_r^{-1} J \ur$, where
$Q_c$ is the coarse ($N=1$) space assembly operator,
$W$ = diag($w_i^{-1}$) is a diagonal weight matrix comprising the inverse
entries of the counting vector, $\uw_L = Q_cQ_c^T \ue_L$, and
$\ue_L = [1 \;\; 1 \; \dots \; 1]^T$ is the local unit vector of length $E n_v$.
(This same weighting procedure is used to enforce continuity at the end
of the Schwarz step, (\ref{eq:asm}) or (\ref{eq:asm1}), and can be combined
in the additive case.)
   By using either of these approaches, one can build $A_r$ through matrix
assembly.  The cell-based localization leads to a significant gain in sparsity,
particularly in 3D.  With the full support for $\Psi$, $A_r$ will have $5^d$
nonzeros per row, as illustrated in Fig. \ref{fig:1D}(a) for $d=1$.  With
either of the compact schemes, the number of nonzeros per row is only $3^d$,
which yields almost a 5$\times$ savings for $d=3$.   Finally, switching from a
width-5 stencil to a width-3 stencil reduces by a factor of two the size of the
separator sets when $A_r$ is permuted by a nested-dissection ordering.  For
Cholesky factorization, this 5-to-3 switch reduces the number of nonzeros by
a factor of four.  For the $XX^T$ factorization discussed below, the number of
nonzeros and the communication overhead, both of which scale with the aggregate
size of the separators, are reduced by a factor of two.

 \section{Results}\label{sec:results}

\begin{figure}[t] \centering
{\setlength{\unitlength}{1.0in}
   \begin{picture}(6.500,2.15)(0,.0)
      \put(0.14,-.03){\includegraphics[height=2.10in]{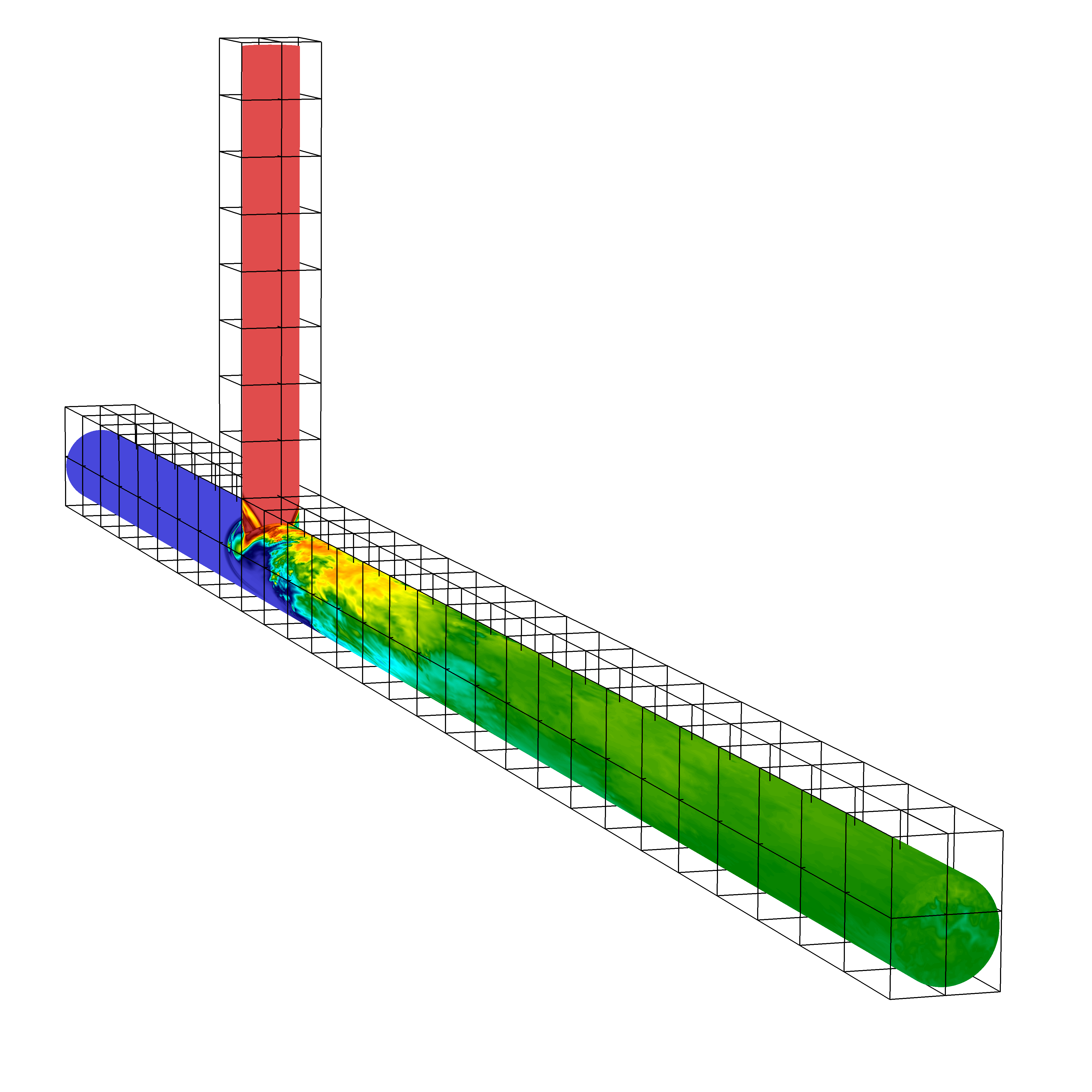}}
      \put(2.10,0.00){\includegraphics[height=2.00in]{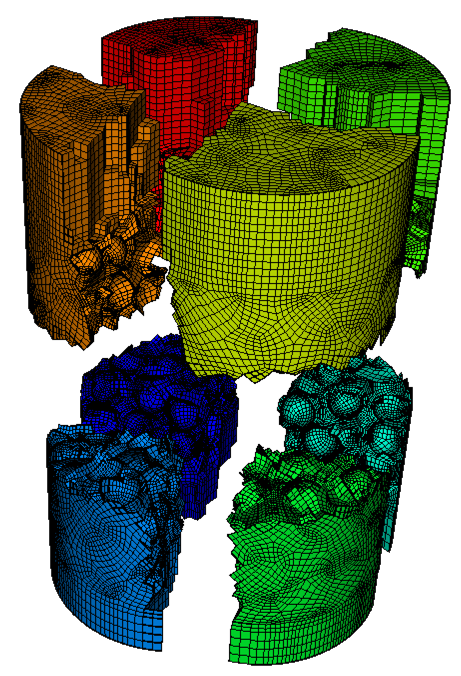}}
      \put(3.40,0.00){\includegraphics[height=2.00in]{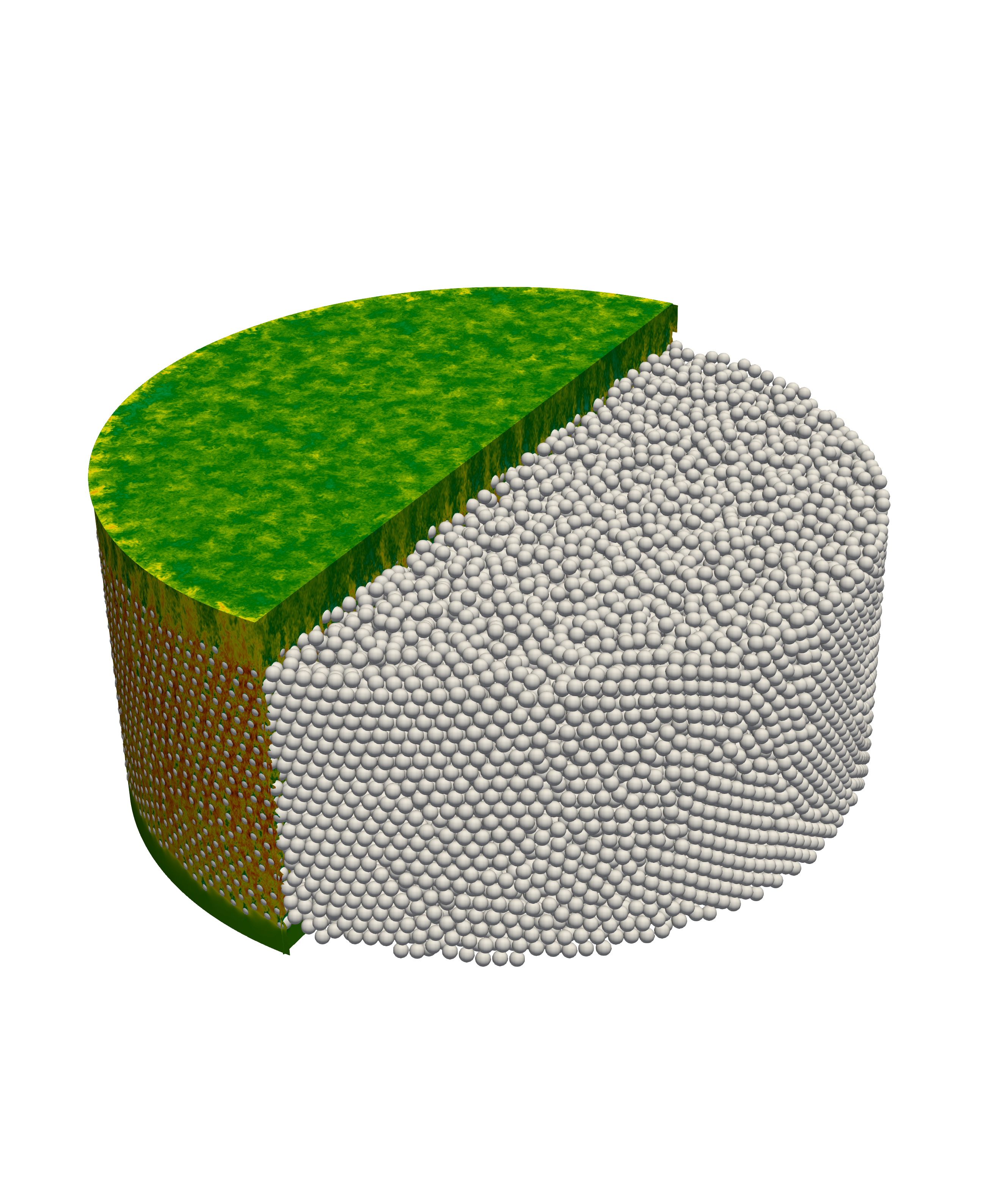}}
      \put(5.00,0.00){\includegraphics[height=1.80in]{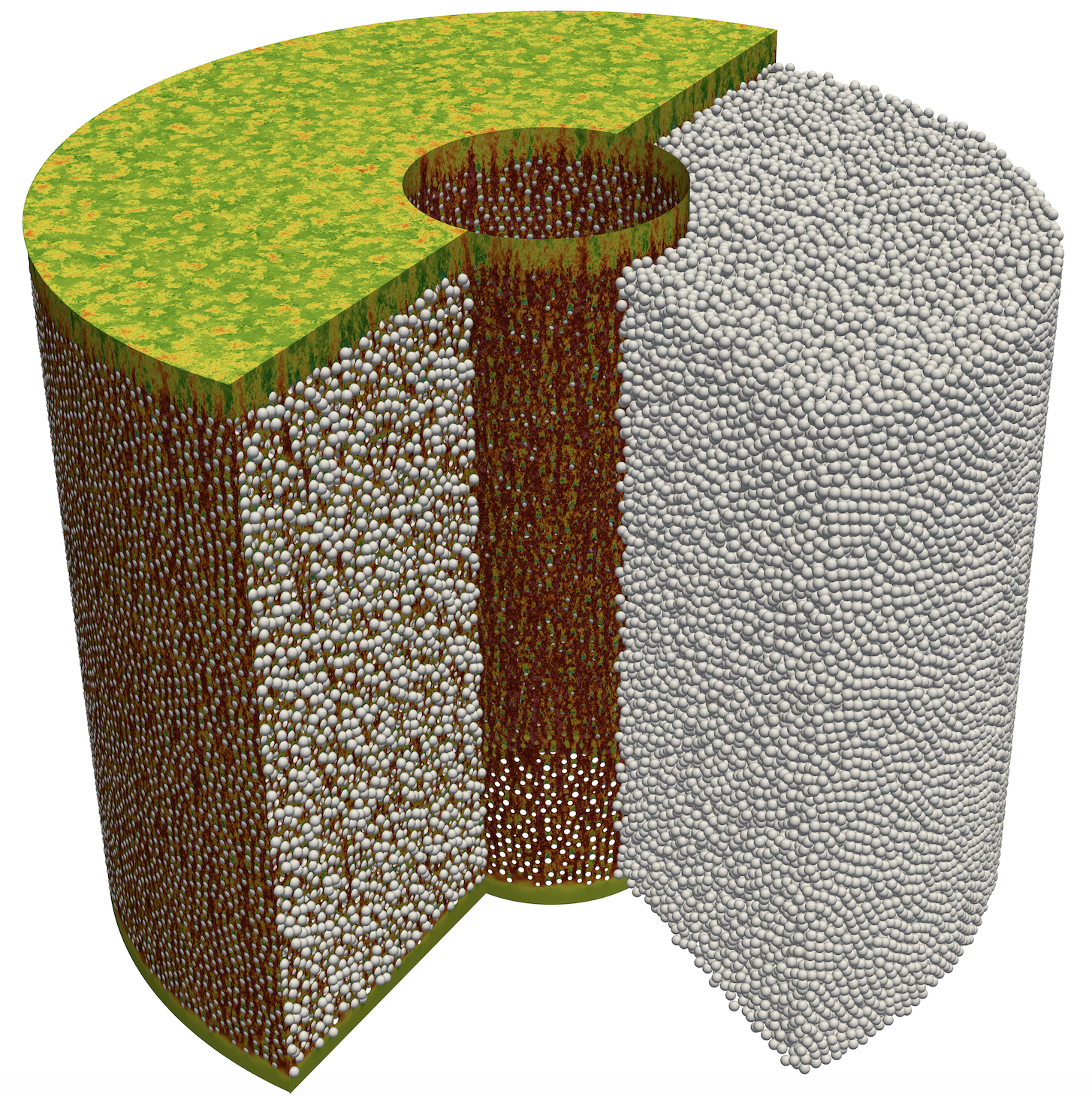}}
   \end{picture}}
   \caption{Test cases (left to right):
            T-junction ($E$=62176) configuration, including coarse-space cells;
            146-pebble ($E$=62138) mesh illustrating an eight-way mesh partition;
            45000-pebble ($E$=13032440) mesh, and 352625-pebble ($E$=98782067) mesh.}
    \label{fig:3dtests}
\end{figure}

In this section, we explore the performance of the new coarse grid solver in the
context of the SEM-based incompressible flow solvers, Nek5000 and NekRS.
A critical question in the current development was whether the outer,
$p$MG-preconditioned, GMRES iteration counts for the pressure solves would change
when the default AMG coarse solver was replaced by the proposed two-level
Schwarz solve. If the reduction of the coarse solve time using the new coarse
grid solver can't compensate for the increased cost due to potentially large outer
iterations, then the new approach will not be competitive.

To answer this question, we compare AMG and the new solver as a coarse grid solver
in four production-level Navier-Stokes simulations.
The example cases used are illustrated in Fig.~\ref{fig:3dtests}.  The first is flow in
a T-junction, which has a mesh comprising $E=62176$ spectral elements of order
$N=7$.  The reduced-space domain is based on a box of having $24 \times 2
\times 7$ cells.  Only a small fraction of these are used, however, as only
those cells that contain spectral element centroids are activated.  These
cells are outlined in the figure.  Note that the majority of the active
reduced-space vertices are exterior to the domain in this case.
The second test case is flow through a 146-pebble bed having $E=62138$ elements
of order $N=7$, which is illustrated in the exploded view in Fig.
\ref{fig:3dtests}. The fluid flows around the spherical pebbles and the
sphere interiors are external to $\Omega$.  Here, we use a $7 \times 7 \times
12$ box mesh for the reduced space, which nominally has 832 DOFs save that
corner cells in the $x$-$y$ plane fall outside the disk, so those corner DOFs
are not included.  This example shows that complex domains with void inclusions
present no particular difficulty for the new coarse-space approach.
The third test is also a pebble bed but with $\approx 45000$ pebbles with
$E=13032440$ elements of order $N=7$. Reduced space is based on a box mesh of
$46 \times 46 \times 23$ for Nek5000 and $28 \times 28 \times 14$ for NekRS.
The forth test is vertical flow through an annular bed of 352625 pebbles
with $E=98782067$ elements of order $N=8$ ($n=51$B).
The reduced space is based on a box mesh of $47 \times 47 \times 35$ for Nek5000
and $32 \times 32 \times 25$ for NekRS.
Again, any box that is fully exterior to $\Omega$ is automatically dropped by the
implementation, which only adds active vertices to the set of reduced-space unknowns.
For each test case, the pressure-Poisson boundary conditions are Neumann
everywhere except at the domain outflow (the end of the long branch for the
T-junction case and the upper boundary for the pebble bed cases), where the
solution has homogeneous Dirichlet conditions.
Although these boundary conditions are not applied to the reduced-space basis
functions, $\Phi_j$, they are automatically imposed on $\Psi_j$.
The original coarse-space dimension is $n_c \approx E$ for all the four cases since
it is the space associated with the last level of $p$MG ($N=1$) which consists of
element vertices.

We run all the tests in the context of full Navier-Stokes simulations,
which require the solution of a pressure Poisson problem
at each time step, $t^m$.  For iterative solvers, it pays to
solve only for the {\em change} in the solution, $\delta \up :=
\up^m - {\bar \up}$, where ${\bar \up} = \up^{m-1}$ or, better, ${\bar \up} =
\sum_{j=1}^l \beta_j \up^{m-j}$ is the $A$-orthognal (best fit) projection onto
$l$ prior solutions \cite{fisc98}.  This projection process tends to remove
low-wavenumber content from the initial residual and typically yields a two- to
three-fold savings in Navier-Stokes run time \cite{fisc98,sc22}.
All Nek5000 runs were performed on the Summit supercomputer at OLCF,
which has 4608 nodes, each with two IBM POWER9 CPUs and six NVIDIA Volta V100 GPUs.
Nek5000 only utilized the CPU cores since it does not support running on the GPUs.
NekRS runs were performed on the Frontier supercomputer at the same facility
which has $\approx 9400$ nodes, each with a single AMD EPYC CPU and four MI250X GPUs.
NekRS was run on the GPUs using HIP backend of OCCA~\cite{occa} portability library
except for the coarse solver which was run on the CPUs in case of two level Schwarz
solver and the default AMG solver using BoomerAMG.

\begin{figure}[t] \centering
{\setlength{\unitlength}{1.0in}
  \begin{picture}(6.500,4.20)(0,0)
   \put(0.0,2.1){\includegraphics[height=2.15in]{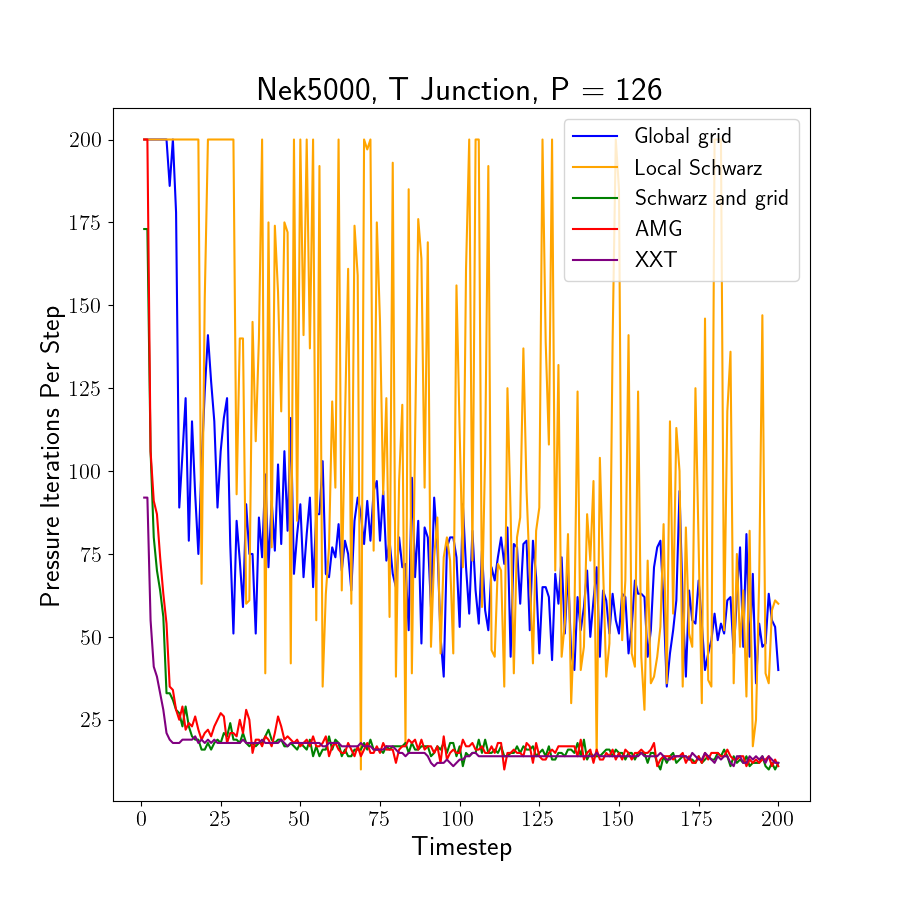}}
   \put(2.2,2.1){\includegraphics[height=2.15in]{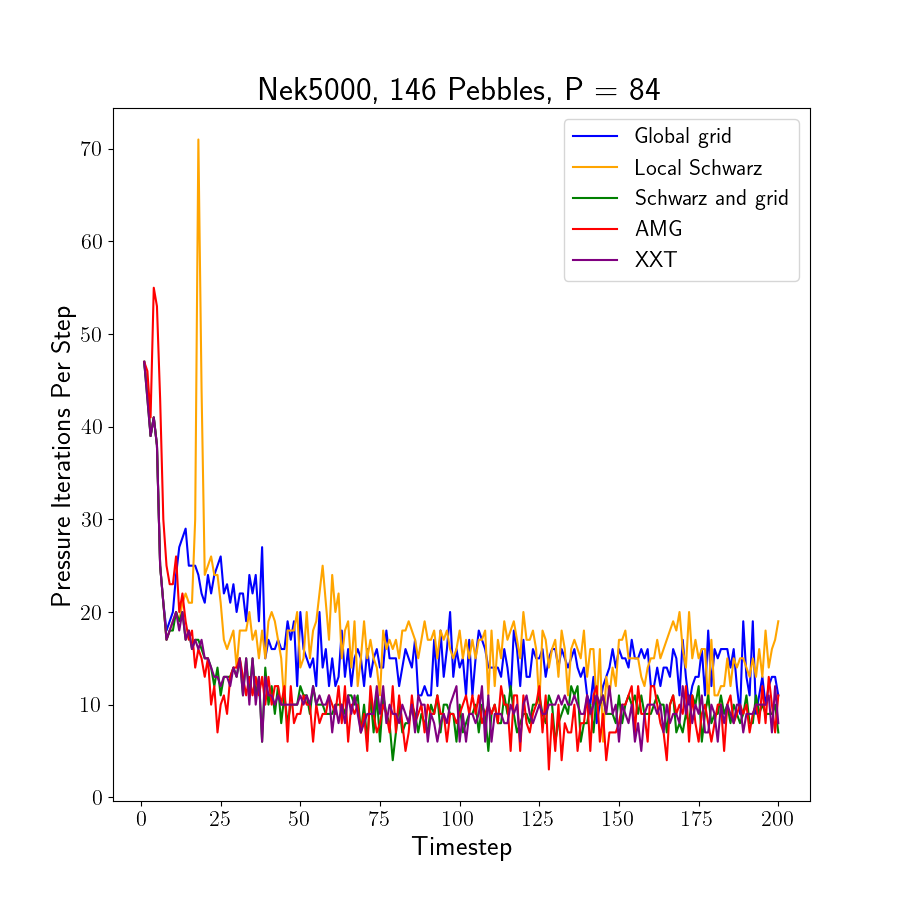}}
   \put(4.4,2.1){\includegraphics[height=2.15in]{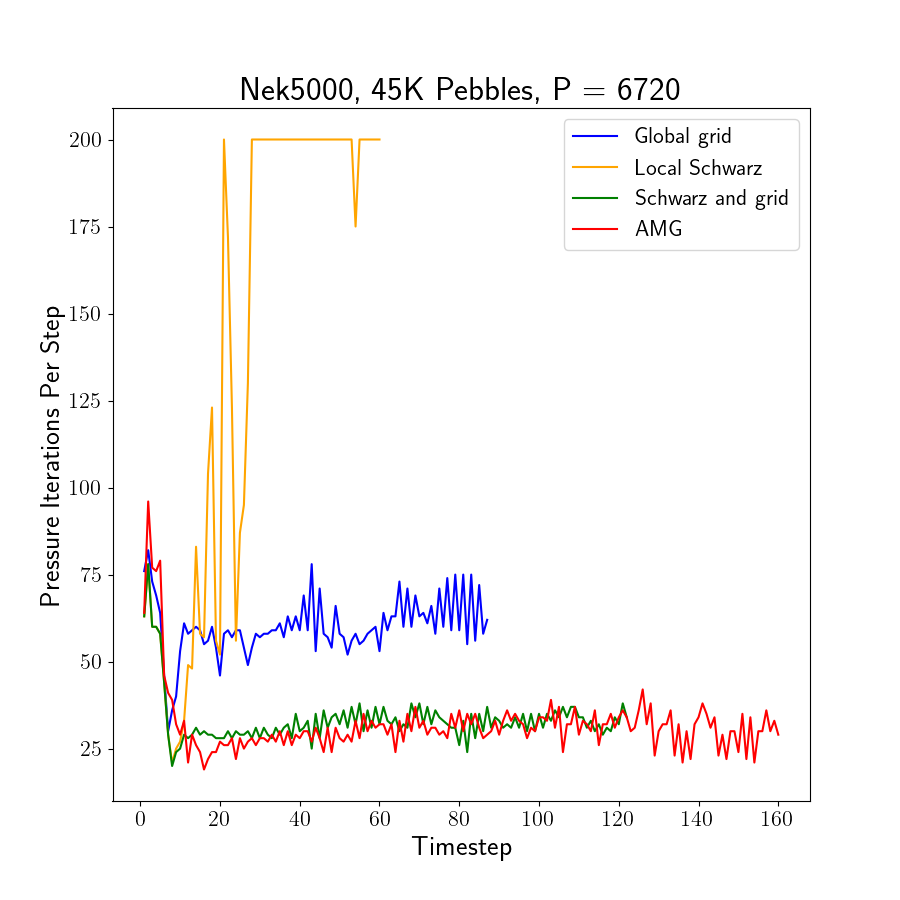}}
   \put(0.0,0.0){\includegraphics[height=2.15in]{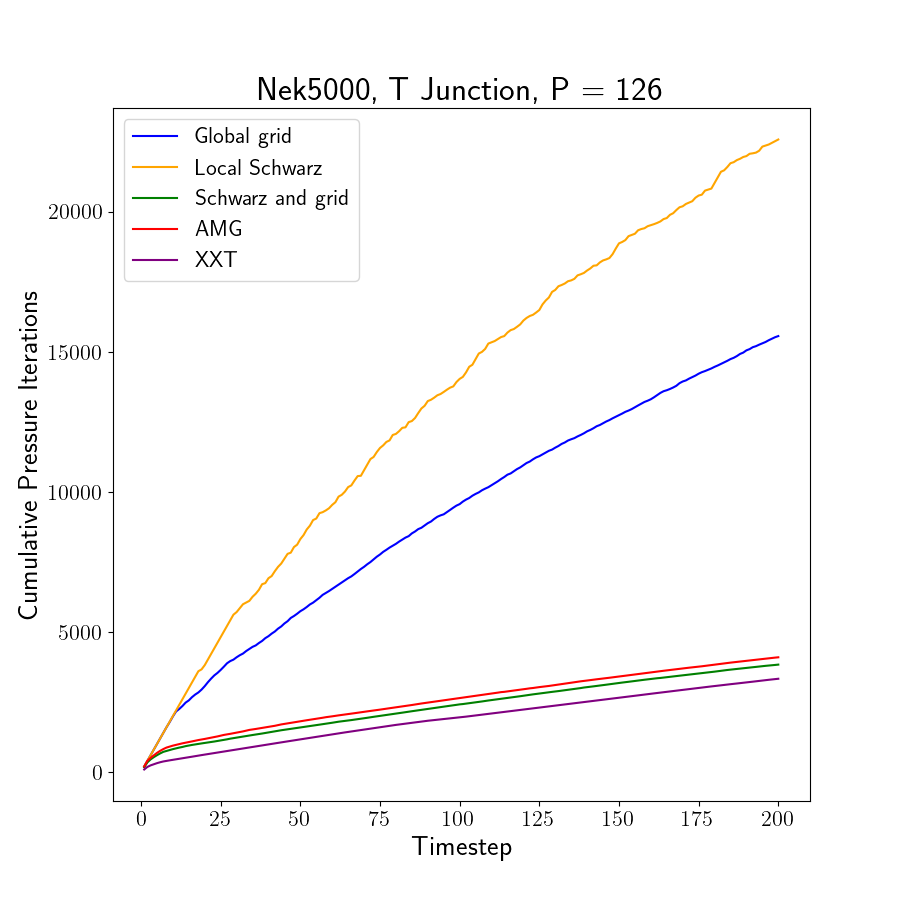}}
   \put(2.2,0.0){\includegraphics[height=2.15in]{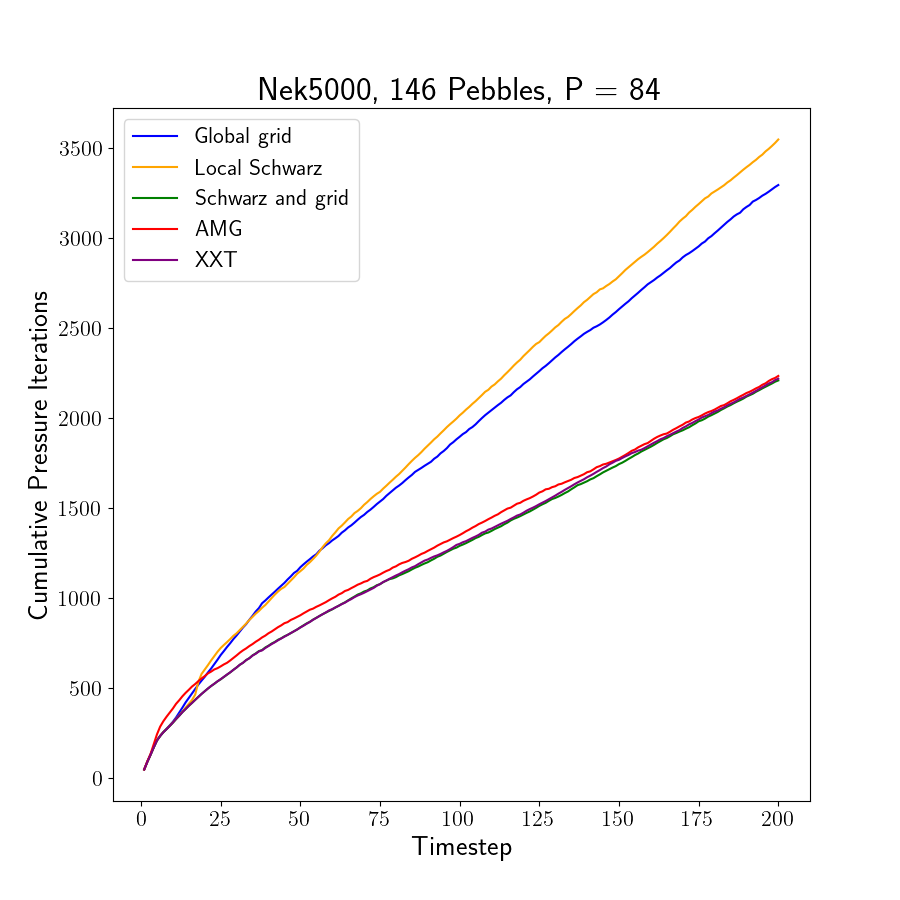}}
   \put(4.4,0.0){\includegraphics[height=2.15in]{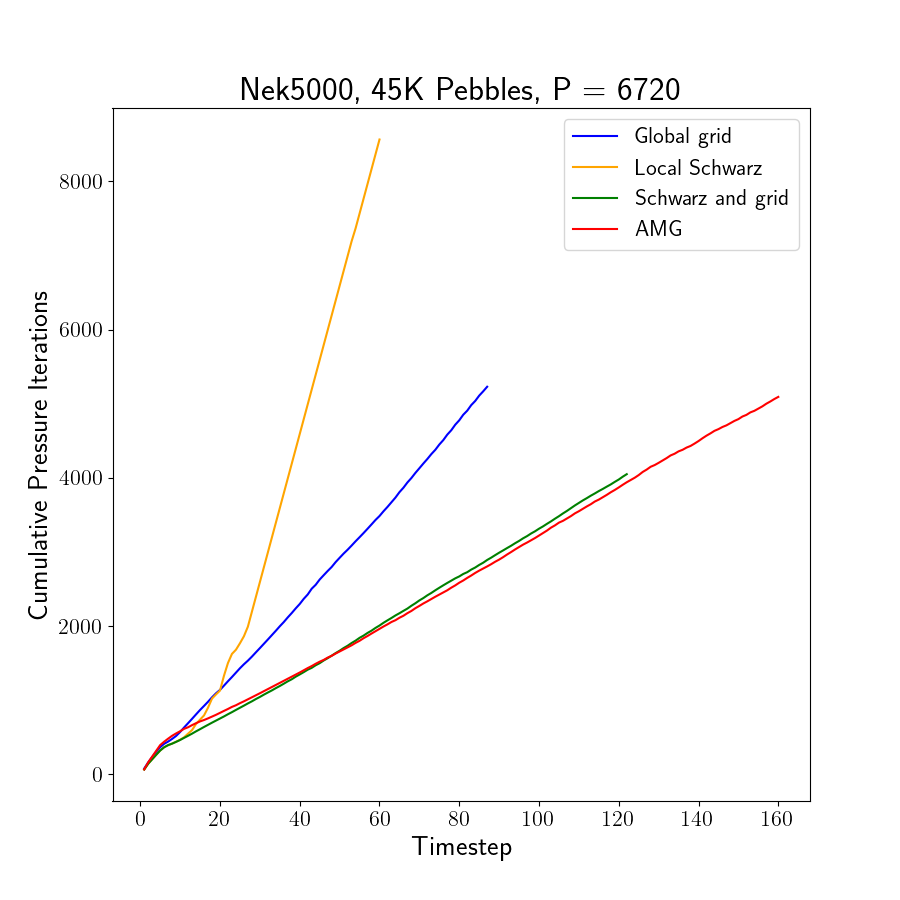}}
  \end{picture}}
  \caption{Nek5000 test results for different coarse solvers:
     (top) iteration count-per-step; (bottom) cumulative iteration counts.
     \label{fig:3diters} }
\end{figure}

The first set of results, shown in Fig.~\ref{fig:3diters}, are for Nek5000,
which runs on CPUs and uses multilevel $p$MG with Schwarz-based smoothing
for the finest two levels, $N=7$ and $N=3$.  The default for the coarse $N=1$
problem is to solve it directly with the $XX^T$ factorization or, for problems
with $E>350000$, to solve it approximately with a single sweep of an AMG
$V$-cycle using Hypre~\cite{boomer-amg-2002}.
The Nek5000
preconditioner, described in \cite{fischer04,lottes05}, is additive between
levels and takes the form
\begin{eqnarray}
  \uz &=& \left(M_N^{-1} \;+\; J_3 \left(M_3^{-1} \;+\; J_1 M_1^{-1} J_1^T \right)J_3^T \, \right) \ur ,
\end{eqnarray}
where $M_1^{-1} =  A^{-1}_c$ represents the coarse solve and
\begin{eqnarray}
  M_N^{-1} &=&  W_N \sum_{e=1}^E R_{N,e}^T \tA_{N,e}^{-1} R_{N,e}
\end{eqnarray}
represents the element-local additive Schwarz smoother for polynomial order $N>1$.
The local Schwarz problems are solved (approximately) using the tensor-product
based fast diagonalization methods described in \cite{fischer00a,lottes05}.  The
diagonal weight matrix, $W_N$, is used to average the solution in the element
overlap regions, which was found to
significantly improve the smoothing properties of the ASM \cite{lottes05}.  A
central idea in the original Nek5000 ASM implementation was to ensure that the
result of each matvec was projected. Consequently, only a single matvec is used
per GMRES iteration, but the overall iteration count is higher than that of
NekRS, which uses multiple smoothing sweeps at each $p$MG level.

Figure~\ref{fig:3diters} shows the pressure iteration counts (or number of coarse
grid solves) per time step for each of the examples when using several different coarse
solvers for the multilevel $p$MG preconditioner in Nek5000. The top row shows the results as
number of iterations per step, whereas the bottom shows the cumulative
iteration counts, which are easier to read given the step-to-step variability
in iteration counts.  (This variability is largely due to the pressure
projection, which extracts all temporal regularity out of the solution such
that the initial residual is highly variable \cite{fisc98} .) Shown are five
different coarse solvers: using just the global reduced-space (box) solver,
$M_r^{-1} := J A_r^{-1} J^T$; using just the local Schwarz solve, ASM$_1$
\eqref{eq:asm1}; combining these into a two-level Schwarz solve, ASM$_2$
\eqref{eq:asm2}; using a single $V$-cycle of Hypre; and using the direct
$XX^T$-based solve.  First off, we note that {\em both parts} of the ASM$_2$
solver are crucial for the success of the method. Neither the local solve nor
the reduced-space box-solve are sufficient to yield a fast $p$MG solver.
Moreover, we see that ASM$_2$, with the sparse box-solver for the coarse grid
problem, is almost as effective as $XX^T$ and slightly more
effective than a single AMG $V$-cycle.
It is remarkable just how effective the low-communication coarse space is for these
applications.

\newpage
\begin{figure}[t] \centering
{\setlength{\unitlength}{1.0in}
  \begin{picture}(6.500,2.10)(0,0)
   \put(0.0,0.0){\includegraphics[height=2.15in]{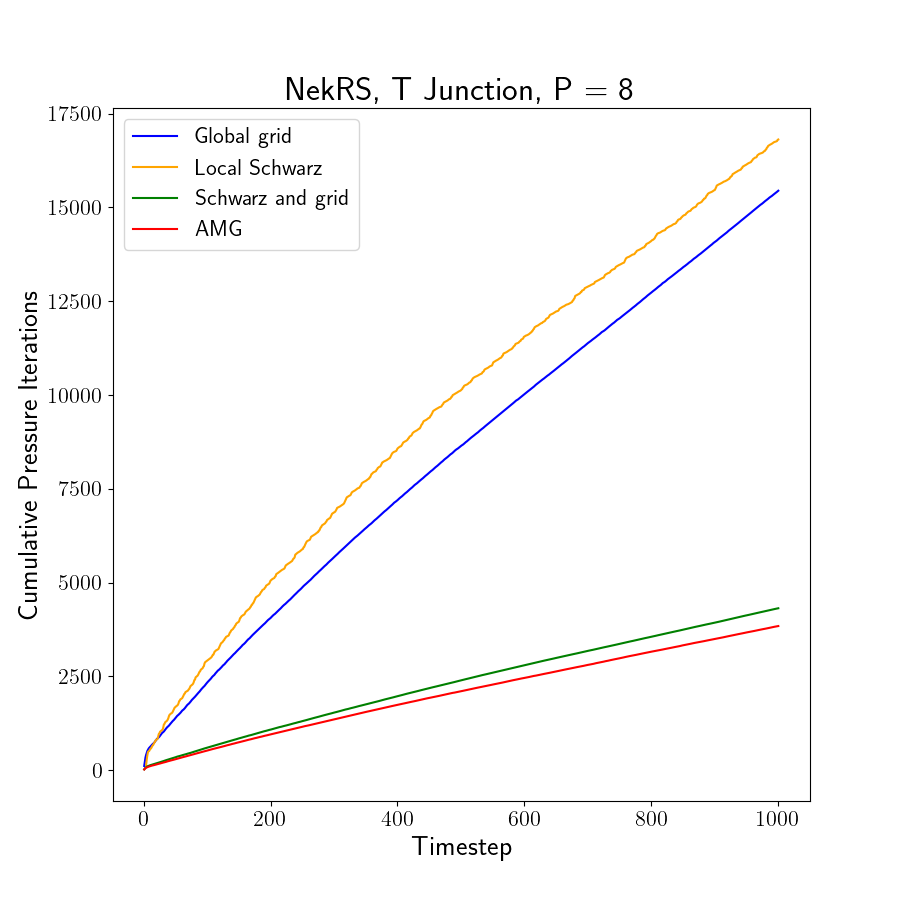}}
   \put(2.2,0.0){\includegraphics[height=2.15in]{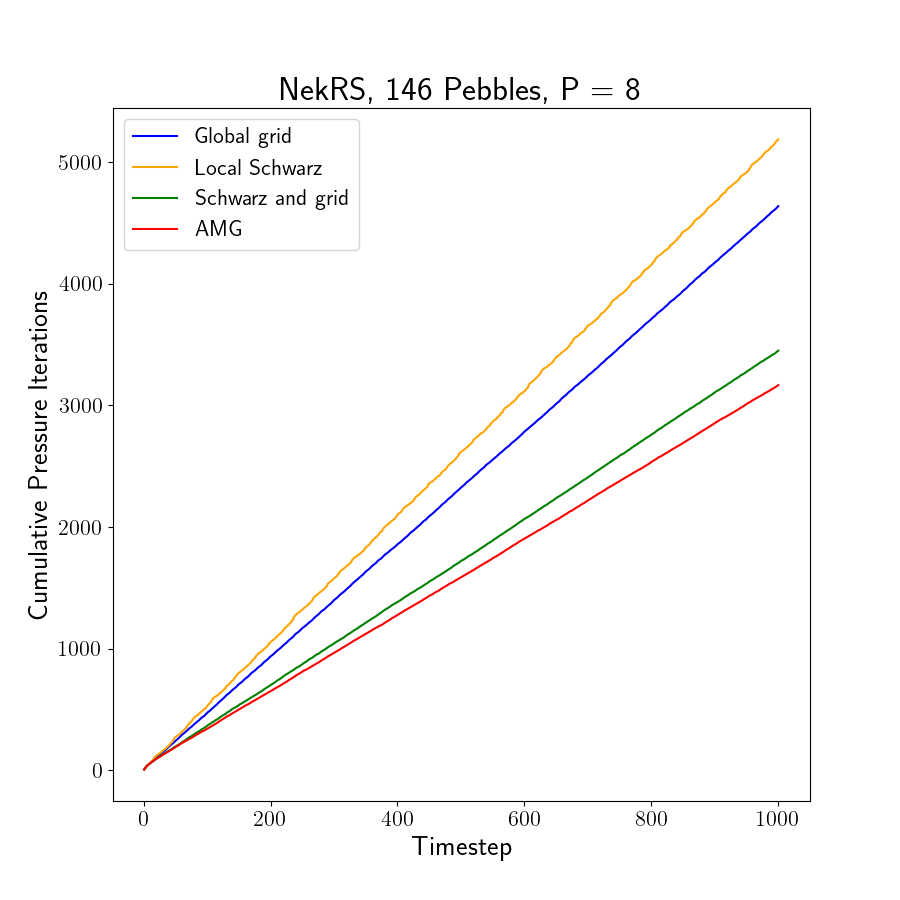}}
   \put(4.4,0.0){\includegraphics[height=2.15in]{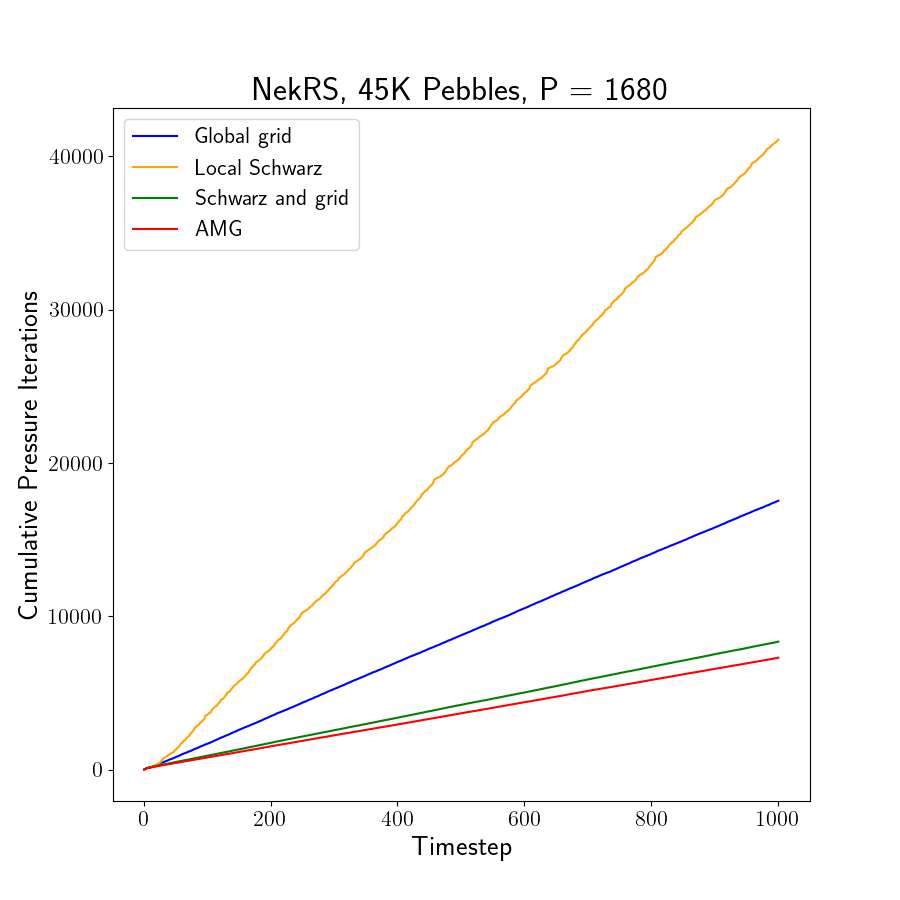}}
  \end{picture}}
  \caption{NekRS results for different coarse solvers: cumulative iteration counts.\label{fig:3diters_nekrs}}
\end{figure}

Second set of results, shown in Figure~\ref{fig:3diters_nekrs} are for the same set of
cases simulated using NekRS running on GPUs except for the coarse solver which is running
on the CPUs.
We only plot the cumulative iterations counts for the cases with NekRS in
Figure~\ref{fig:3diters_nekrs}.
NekRS uses a multilevel $p$MG with a Schwarz based smoother and a schedule $N=7\to3\to1$
similar to Nek5000 settings used in experiments in Figure~\ref{fig:3diters}.

Default coarse solver in NekRS is a single sweep of an AMG $V$-cycle on the original
system $A_c$ ($N=1$) using Hypre~\cite{boomer-amg-2002} in single precision on the CPU.
Similar to what we saw with Nek5000, both parts of ASM$_2$ are required for the
success of the new coarse solver with NekRS as well.
One main difference compared to Nek5000 experiments is that a single $V$-cycle with Hypre
is more effective in terms of the cumulative pressure iteration counts than the Schwarz
based coarse solver for the NekRS experiments.
Our primary interest is the time per timestep at large scale simulations rather than the
cumulative iteration counts and the whole point of the design of the new Schwarz based
coarse solver was to have a smaller cumulative coarse grid solve time compared to AMG
albeit a slightly larger cumulative number of pressure iterations.

\begin{table}
  \centering
  \caption{Comparison of solver times and pressure iterations per timestep with NekRS when using
  AMG and two level Schwarz solver as a coarse grid solve.}
  \begin{tabular}{|l|l|r|r|r|r|r|r|r|}
    \hline
    Case & Solver & P & E/P & n/P & NS Time(s) & Pres. Time(s) & Coarse Time(s) & Pres. Iter. \\
    \hline
    T-Junction  & Schwarz&     8& 7772 & 2.66e+06 & 1.37e-01 & 8.59e-02 & 2.57e-02 & 4.32 \\
    \hline
    T-Junction  & AMG    &     8& 7772 & 2.66e+06 & 1.21e-01 & 6.97e-02 & 1.78e-02 & 3.84 \\
    \hline
    146 Pebbles & Schwarz&     8& 7767 & 2.66e+06 & 1.14e-01 & 6.72e-02 & 1.53e-02 & 3.45 \\
		\hline
    146 Pebbles & AMG    &     8& 7767 & 2.66e+06 & 1.14e-01 & 6.63e-02 & 1.68e-02 & 3.17  \\
    \hline
    45K Pebbles & Schwarz&  1680& 7758 & 2.66e+06 & 2.60e-01 & 1.95e-01 & 4.66e-02 & 8.35 \\
    \hline
    45K Pebbles & AMG    &  1680& 7758 & 2.66e+06 & 2.69e-01 & 2.05e-01 & 7.35e-02 & 7.31 \\
    \hline
    350K Pebbles& Schwarz& 27648& 3573 & 1.83e+06 & 2.17e-01 & 1.59e-01 & 1.70e-02 & 6.18 \\
    \hline
    350K Pebbles& AMG    & 27648& 3573 & 1.83e+06 & 2.18e-01 & 1.61e-01 & 4.56e-02 & 4.80 \\
    \hline
  \end{tabular}
  \label{tab:solver_performance}
\end{table}

Table~\ref{tab:solver_performance} compare the solve times and pressure iterations per
timestep when using AMG and two level Schwarz method as the coarse grid solver in NekRS
for four cases including three cases used in experiments in
Figure~\ref{fig:3diters}-\ref{fig:3diters_nekrs} and the annular bed of 352625 pebbles.
The latter case uses a schedule $N=8\to6\to4\to1$ for the multilevel $p$MG preconditioner.
The results show that the two level Schwarz method is competitive
(and does slightly better) at higher process counts where the communication cost of AMG
becomes significant due to multiple levels in the $V$-cycle.
This can be seen for the two larger meshes run with $P=1680$ and $P=27648$ in Table~\ref{tab:solver_performance}
where the time spent in coarse grid solve per timestep is about 1.6 and 2.7 times faster
than default AMG solver respectively.
These speedups eventually translate to an overall speedup of full Navier-Stokes solve time
per timestep.
While, the number of pressure iterations per timestep is lower with AMG for all the cases,
reduced coarse times eventually helped the two level method catch up with AMG at higher
process counts for the two large cases with high number of MPI processes.

\begin{figure}[t] \centering
{\setlength{\unitlength}{1.0in}
  \begin{picture}(6.500,2.10)(0,0)
   \put(0.0,0.0){\includegraphics[height=2.15in]{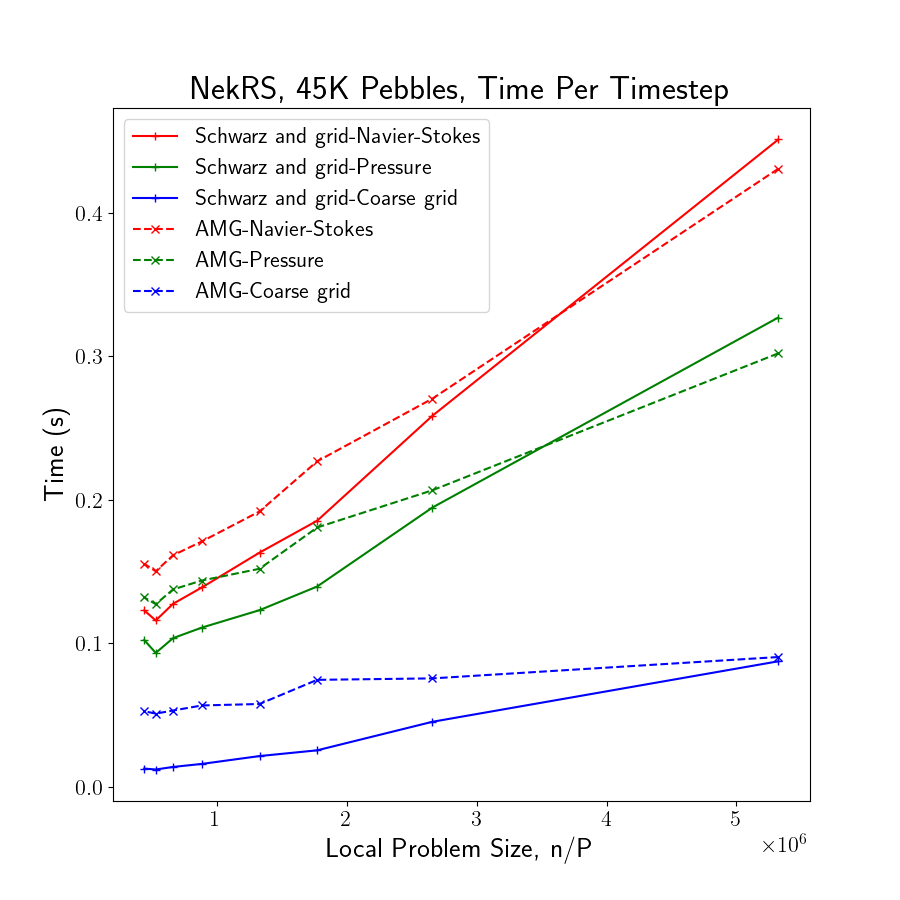}}
   \put(2.2,0.0){\includegraphics[height=2.15in]{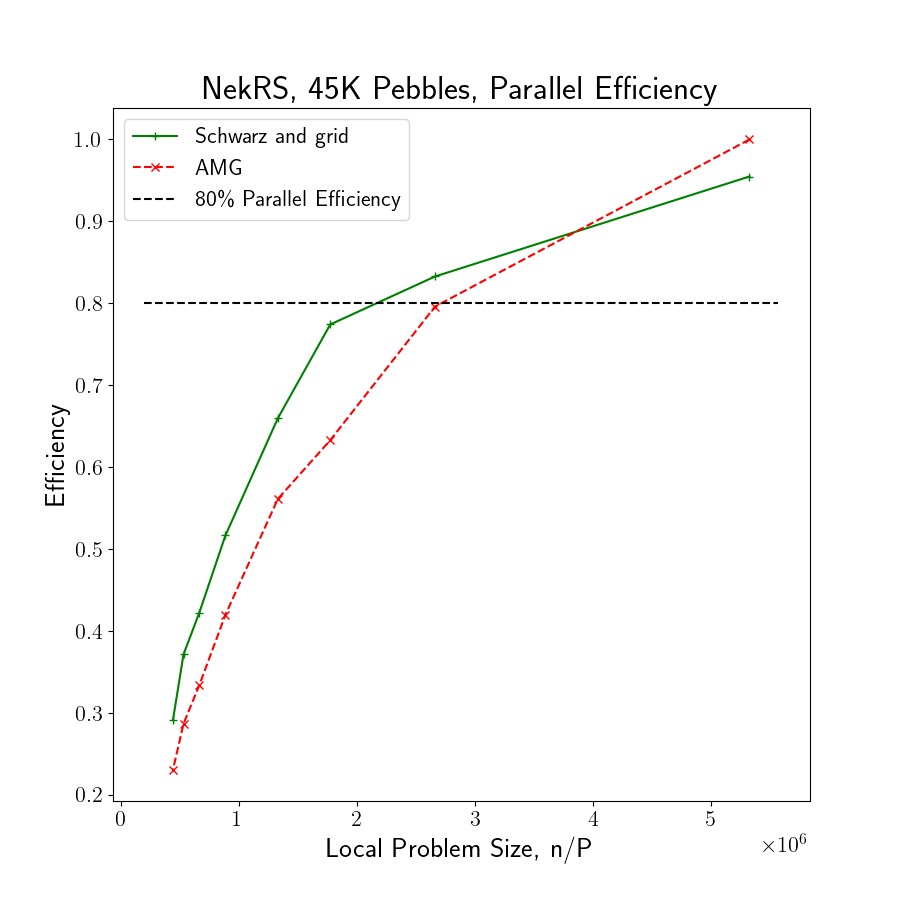}}
   \put(4.4,0.0){\includegraphics[height=2.15in]{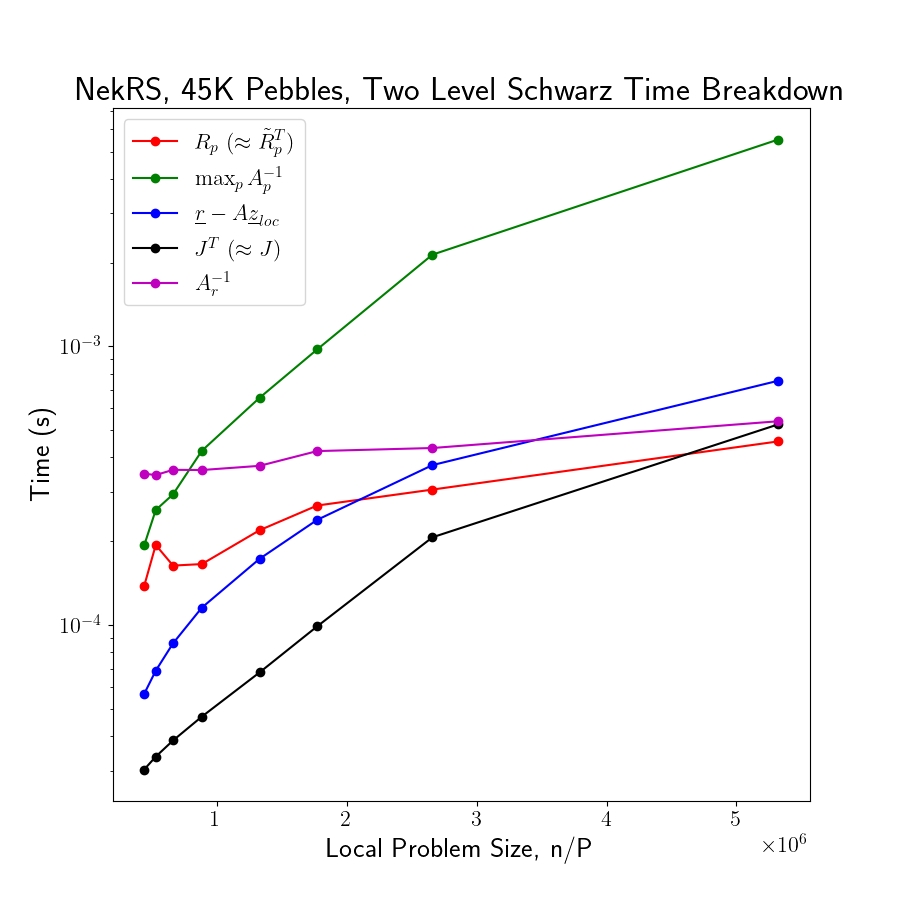}}
  \end{picture}}
  \caption{Strong scaling study of 45000 pebbles mesh with NekRS. Left to right: Navier-Stokes,
  Pressure and Coarse grid solve time per timestep for AMG and two level Schwarz method; Parallel
  efficiency for AMG and two level Schwarz method; Breakdown of the coarse grid solve time for the
  two level Schwarz method.}
  \label{fig:45k_peb_ss}
\end{figure}


Figure~\ref{fig:45k_peb_ss}, which shows timing data for a strong scaling study of the 45000
pebbles mesh with NekRS further demonstrate this point.
The plots were generated by collecting timing data by increasing the number of processes, $P$,
from $840$ to $10080$ (or decreasing $n/P$).
We can see that Navier-Stokes solve time, pressure solve time and coarse grid solve time per
timestep goes down faster with the two level Schwarz solver compared to AMG as the number of
processes increase in the leftmost plot in Figure~\ref{fig:45k_peb_ss}.
Crossover point where the two level Schwarz solver becomes faster than the default AMG approach
is $n/P \approx 3.5\times10^6$.
Middle plot of Figure~\ref{fig:45k_peb_ss} shows the parallel efficiency for the two solvers as
$n/P$ decreases.
Reference point for the parallel efficiency is the NekRS run with Hypre AMG solver with $P=840$
which is assumed to have 100\% parallel efficiency.
It is evident from the middle plot that the new two level Schwarz solver has better parallel
efficiency than AMG for $n/P$ values which falls in our original design range for $n/P$.

Right most plot in Figure~\ref{fig:45k_peb_ss} shows the breakdown of the coarse grid solve
time for the two level Schwarz solver.
Cost of the operations which don't involve communication goes down noticeably as $n/P$ decreases.
These include the local Schwarz solve $(A_p^{-1})$, right hand side update
$(\underline{r} - A \underline{z}_{loc})$ required for doing reduced solve in a multiplicative
manner and interpolation $(J)$ from the original coarse space $A_c$ to the reduced structued
space $A_r$.
Note that in the case of $J$ and $J^T$, it is the structured nature of $A_r$ which allows us
to compute the interpolation and prolongation in parallel with no communication.
The fact that the local Schwarz solve is the most expensive operation in the two level Schwarz
for most of the range presents an opportunity for further optimizations which are discussed breifly
in Section~\ref{sec:future}.
Cost of solving the reduced system $A_r^{-1}$ stays more or less the same as $n/P$ decreases.
This is not suprising since the size of this system is not affected by $n/P$ but fixed for the
entire study.
Althoug a slight decrease in the cost of solving $A_r^{-1}$ is observed as $n/P$ decreases
which hints at the fact that it is the local part of $XX^T$ which is dominating the cost of
solving $A_r^{-1}$.

Table~\ref{tab:350k_chebyshev} illustrate the importance of visiting the coarse grid
system as few times as possible in order to reduce the overall Navier-Stokes solve time.
This is the strategy adapted in~\cite{sc22} in order to reduce the coarse grid overhead
indirectly by doing more Chebyshev smoothing steps in order to reduce the number of
pressure solves. Increasing number of Chebyshev smoothing steps from 1 to 3 decreases the
number of pressure iterations by more than a factor of three.

\begin{table}
  \centering
  \caption{Navier-Stokes solve times per timestep with different Chebyshev orders for 350K pebbles
  mesh with NekRS using $P=27648$ processes with $\approx 2M$ grid points per process when using
  AMG and two level Schwarz solver as coarse grid solvers.}
  \begin{tabular}{|l|r|r|r|r|r|}
    \hline
     Solver & Cheb. Order & NS Time(s) &   Pres.Time(s) &   Coarse Time(s) &   Pres. Iter. \\
     \hline
     AMG    &           3 &   2.18e-01 &       1.61e-01 &         4.56e-02 &      4.80 \\
     \hline
     AMG    &           2 &   2.60e-01 &       2.02e-01 &         7.02e-02 &      7.28 \\
     \hline
     AMG    &           1 &   3.65e-01 &       3.04e-01 &         1.38e-01 &      14.50 \\
     \hline
     Schwarz&           3 &   2.17e-01 &       1.59e-01 &         1.70e-02 &      6.18 \\
     \hline
     Schwarz&           2 &   2.48e-01 &       1.89e-01 &         2.77e-02 &      9.78 \\
     \hline
     Schwarz&           1 &   3.31e-01 &       2.75e-01 &         5.68e-02 &      20.90 \\
    \hline
  \end{tabular}
  \label{tab:350k_chebyshev}
\end{table}

 \section{Discussion}\label{sec:discussion}

In this section, first we are going to do a complexity analysis of both AMG and
two level Schwarz method to see if we are on the right ballpark in terms of the
timing data collected and to estimate the costs for future exascale runs.
Then we are going to analyse the affect of the new two level Schwarz coarse solver
on strong scaling of Nek5000/RS.
Finally, we are going to discuss the future work that we are planning to do in
order to improve the performance of the solver.

\subsection{Complexity Analysis}\label{sec:complexity}
\begin{center}
\begin{table}
\begin{tabular}{|c|c|c|}
\hline
 Description & Symbol & Measured or Estimated value\\ \hline
 Latency of the Network& $\alpha$ & 4e-6 \\ \hline
 Short-long message demarcation & $m_2 = \alpha/\beta$ & $\approx 5000$ \\ \hline
 Arithmetic time (inverse FLOPS) & $t_a$ & $10^{-9}\mbox{s}$ \\
\hline
\end{tabular}
\caption{System parameters for Frontier supercomputer at OLCF}
\label{tab:sys_param}
\end{table}
\end{center}

Here we develop time complexity estimates for parallel evaluation of the two-level
Schwarz and AMG in the context of a coarse grid solver for pMG preconditioner.
We are considering a problem size consitent with the strong scaling limit for GPUs,
which is about $n/P \approx 2.66\times10^6$ or $E/P \approx 8000$  with $N=7$.
For Schwarz, we assume that the local and reduced coarse-space solves are evaluated
serially in a multiplicative way, so that one might consider using a hybrid Schwarz
approach in which the additive local solves are followed by a coarse grid correction.
We define the required system parameters, along with estimated/measured values for
Frontier supercomputer at OLCF in Table~\ref{tab:sys_param}.
Note that both AMG and Schwarz are run on the CPU in the context of the coarse
solver, so we only need the parameters for the CPU and the network.

For ASM$_1$, we have pre- and post-solve ($R_p$ and $\tilde{R}^T$) near-neighbor
communication requiring $\approx$ 50 messages of size $\leq (E/P)^{2/3}$.
For GPUs at the strong-scale limit, we anticipate $E/P \approx 8000$ and therefore
message sizes $m \leq (E/P)^{2/3} = 400 < m_2$, which implies that the messages are
latency bound, each with a cost $\leq 2 \alpha$.
Thus, anticipated communication cost for ASM$_1$ is
\begin{eqnarray}
  T_{c,ASM_1} &=& 50 \cdot (2 \alpha) = 100 \cdot \alpha  \approx .0002\mbox{s}
   \label{eq:asm1_comm}
\end{eqnarray}
Above communication time is the sum of the time for $R_p$ and $\tilde{R}^T$ shown
in the rightmost plot of Figure~\ref{fig:45k_peb_ss}.
Cost of $R_p$ and $\tilde{R}^T$ are more or less the same, and thus should be half of
the above estimate.
We can see that the timing data for $R_p$ and $\tilde{R}^T$ for $E/P \approx 8000$
in rightmost Figure~\ref{fig:45k_peb_ss} is in the same order of magnitude as the
estimate above.
Also, there is a local computation to be done for $R_p$ and $\tilde{R}_p^T$ in addition
to the communication which explains why the measured time goes down slightly as local
problem size $n/P$ goes down.

Local solve for ASM$_1$ is implemented as a Cholesky solve on the CPU, so we can
estimate the time using the CPU arithmetic time $t_a$ and the number of operations
performed during the for the local solve.
Assuming a Cholesky factorization of $A_p = LL^T$, the number of operations in forward
and backward solve is roughly $2n_z$ where $n_z$ is the number of nonzeros in $L$
(we have to access each non-zero value in $L$ or $L^T$ once during the solve and
update the right hand side by substracting its contribution after the multiplication by
relevant solution component).
We have noticed that $n_z$ could be as high as  $\num{5e6}$ for the local problems when
$E/P\approx8000$ (or $n/P \approx 2.66\times10^6$ for $N=7$).
\begin{eqnarray}
  T_{a,ASM_1} &=& 2n_z t_a \approx 0.01\mbox{s}
  \label{eq:asm1_arith}
\end{eqnarray}

With both communication in Equation~\ref{eq:asm1_comm} and arithmetic time in
Equation~\ref{eq:asm1_arith}, we can estimate the total time for ASM$_1$ as:
\begin{eqnarray}
  T_{ASM_1} &=& T_{c,ASM_1} + T_{a,ASM_1} \approx 0.01\mbox{s} \label{eq:asm1_total}
\end{eqnarray}

For the reduced-space solve $A_ru_r=b_r$, we follow the complexity estimate for the
$XX^T$ solver developed in \cite{tufofisc01}. For a regular 3D mesh (which provides a
conservative bound), the estimated communication and arithmetic times for $XX^T$ are:
\begin{eqnarray}
  T_{c, A_r^{-1}} &=& 2 \left[\left(1+\frac{2.5}{m_2}n_r^{2/3}\right)\log_2 P\right]\alpha \label{eq:xxt_comm}\\
  T_{a, A_r^{-1}} &=& 2 \, \left( 2.5 \, n_r^{5/3}/P \right) \, t_a \label{eq:xxt_arith}
\end{eqnarray}
For the 45000 Pebbles problem, reduced system size $n_r$ is $28\times28\times14\approx11000$
for NekRS.
If we consider $P=1680$ (for which $E/P\approx 8000$ and $n/P=2.66\times10^6$), we can
calculate the arithmetic, communication and total time for the reduced-space solve as follows:
\begin{eqnarray}
  T_{c,A_r^{-1}} &=& 2\left[\left(1 + \frac{2.5}{5000}\times11000^{2/3}\right)\times\log_2 1680\right]\times\num{4e-6}\\
  &=&2\left[\left(1+\frac{2.5}{5000}\times500\right)\times11\right]\times\num{4e-6}=2\times13.5\times\num{4e-6}=0.0001 \mbox{s}\\
  T_{a,A_r^{-1}} &=& 2 \, \left(2.5\times11000^{5/3}/1680\right)\times\num{1e-9}
  =2\times\left(2.5\times5.4e6/1680\right)\times\num{1e-9}\\
  &=&2\times\left(2.5\times3200\right)\times\num{1e-9}=0.00001\mbox{s}\\
  T_{A_r^{-1}} &=& T_{c,A_r^{-1}} + T_{a,A_r^{-1}} = 0.0001\mbox{s}
\end{eqnarray}
From leftmost plot in Figure~\ref{fig:45k_peb_ss}, we can see that the measured time for
the reduced-space is about $0.0003$s which is the same order of magnitude as the complexity
estimate.
Also, we expect the actual run time to be higher than the estimate since we don't have a
nested dissection ordering for the reduced system and the Equation~\ref{eq:xxt_arith} is
derived assuming a nested dissection ordering.

We contrast these times with estimated and measured times for AMG applied to the full
system $A\uu = \ub$ (i.e., $A_c \uu_c = \ub_c$ in the pMG solver).
BoomerAMG based coarse solver used in the experiments performed a single $V$-cycle with
a single smoothing step in each level using a Chebyshev smoother with degree $2$.
Thus, each smoothing step consists of two nearest neighbor communications with about $26$
messages per communication.
These messages are also latency bound similar to $R_p$ and $\tilde{R}_p^T$ in the ASM$_1$
solver with a cost of $2\alpha$.
If there are $L$ levels in the $V-$cycle, the communication cost of the full $V-$cycle
is given by:
\begin{eqnarray}
  T_{c,AMG} &=& 2 \cdot L \cdot (2 \cdot 2 \alpha \times 26) \approx 200\cdot L\cdot \alpha\label{eq:amg_comm}
\end{eqnarray}
The total arithmetic time is estimated as twice the cost of the smoothing step on the
fine grid which is a sparse matvec with $E/P$ rows and $\approx 27$ non-zeros per row.
\begin{eqnarray}
  T_{a,AMG} &=& 2 \cdot \left( 27 E/P \right) \cdot t_a \approx 50 \cdot E/P \cdot t_a \label{eq:amg_arith}
\end{eqnarray}
Thus, the total time for AMG solver with a sing $V-$cycle is given by:
\begin{eqnarray}
  T_{AMG} &=& T_{c,AMG} + T_{a,AMG} = 200\cdot L\cdot \alpha + 50\cdot E/P\cdot t_a \label{eq:amg_total}
\end{eqnarray}

For the 45000 Pebbles problem, at $P=1680$, BoomerAMG had $L=10$ levels in the
$V-$cycle and the local problem $E/P \approx 8000$ similar to the Schwarz based
approach.
Communication and arithmetic times for the $V-$cycle can then be calculated as:
\begin{eqnarray}
  T_{c,AMG} &=& 200\cdot10\cdot\alpha = 0.008\mbox{s}\\
  T_{a,AMG} &=& 50\cdot8000\cdot\num{1e-9} = 0.0004\mbox{s}\\
  T_{AMG} &=& T_{c,AMG} + T_{a,AMG} = 0.0084\mbox{s}
\end{eqnarray}
Using the data in Table~\ref{tab:solver_performance} for 45K Pebbles problem with AMG,
we can calculate the time for a single coarse solve by dividing coarse time by the number
of pressure iterations which comes out to $7.35e-2/7.31 \approx 0.01$s.
We can see that the measured time for AMG matches pretty well with the estimated time.

\subsection{Strong Scaling}\label{sec:ss}

In this section, we investigate strong scalability of Schwarz and AMG based
coarse-grid solvers for pMG both experimentally and using the timing estimates
developed in Section~\ref{sec:results} and show that the Schwarz has better
strong scaling properties than AMG.

As shown in Section~\ref{sec:results}, strong scaling experiments done on
45000 pebbles case shown in Figure~\ref{fig:45k_peb_ss} show that, as we
increae the number of processes $P$ (thus decreasing $n/P$), time spent in
both Schwarz and AMG (leftmost plot in Figure~\ref{fig:45k_peb_ss}) decreases
with Schwarz time decreasing faster than AMG.
Middle plot in Figure~\ref{fig:45k_peb_ss} shows the parallel efficiency for the
same experiment assuming BoomerAMG run with $P=840$ has 100\% of effiiciency.
An HPC user running production level simulations has to decide on the trade-off
between parallel efficiency and the time to solution.
It is clear from the plot that as we keep adding processors, time to solution
keep decreasing till the parallel efficiency reaches 30-40\%.
But a low level of parallel efficiency may not be acceptable due to underutilization
of computing resources.
Let's assume that users are willing to accept a parallel efficiency of 80\%.
Then, from the middle plot in Figure~\ref{fig:45k_peb_ss}, we can see that we
can use larger number of processors with Schwarz to achieve a smaller
time to solution than AMG.
We can see from the middle plot in Figure~\ref{fig:45k_peb_ss} that AMG needs
about \num{3e6} gridpoints for 80\% parallel efficiency whereas Schwarz needs
$\approx \num{2e6}$.

\begin{figure}[htbp]
  \centering
  \includegraphics[width=0.5\textwidth]{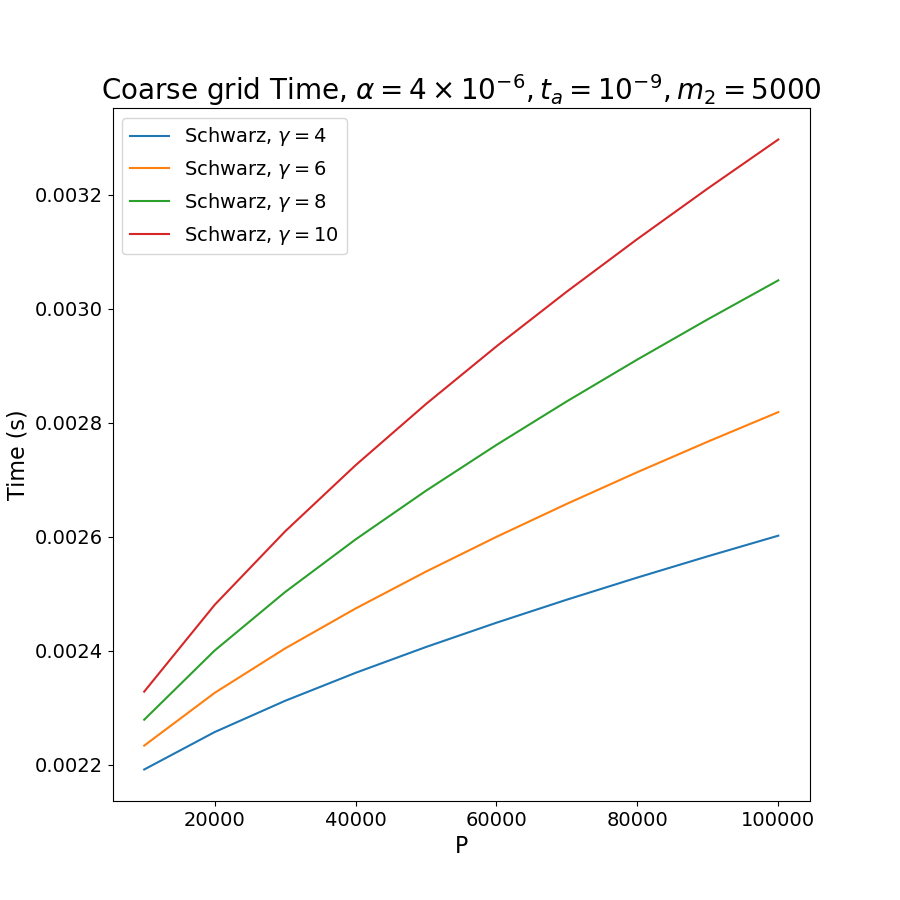}
  \caption{Estimated solve times for Schwarz and AMG for different $\gamma$ and $L$ values.}
  \label{fig:time_estimates}
\end{figure}

We can verify the above observation by looking at the timing estimates developed
in Section~\ref{sec:results}.
Note that we are operating in the regime where local problem size, $E/P\approx 8000$
and $P \approx 10^5$.
For estimating time with AMG, we can use Equation~\ref{eq:amg_total} with $L\approx 10$
since we expect the number of levels in AMG $V$-cycle to be around 10 for problems of
size $E=10^9$.
For the Schwarz solver, our local solve times are going to stay constant since
both local problem size and number of neighbors more or less stays constant as
the number of processors grow.
Only the reduced solve time will increase as we increase the number of processors.
We can calculate the time spent in the reduced solve as the sum of Equation~\ref{eq:xxt_comm}
and Equation~\ref{eq:xxt_arith}.
The reduced solve time is dependent on the number of reduced-space points $n_r$.
By construction, $n_r = \gamma P$, where $\gamma=3$--$10$ is the target number of
reduced-space points per rank.
We assume $\alpha$, and $t_a$ values from Table~\ref{tab:sys_param} for Frontier at
OLCF will also hold true for future exascale systems.
Figure~\ref{fig:time_estimates} show the estimated solve times for Schwarz and AMG
for different $\gamma$ and $L$ values.

\subsection{Conclusion}\label{sec:future}

We presented a two level Schwarz method as a coarse solver for the
$p$-multigrid preconditioner in SEM/FEM simulations with a novel nonnested
coarse space for the global coarse grid of the Schwarz method.
Structured nature of this coarse space enables communication-free interpolation
between the local problem and global coarse grid.
We used CHOLMOD~\cite{chen2008algorithm}, a sparse Cholesky solver for
solving the local problem exactly.
CHOLMOD uses Approximate Minimum Degree (AMD)~\cite{amestoy2004algorithm}
ordering to reduce fill-in in the Cholesky factor.
Global coarse problem is solved exactly using $XX^T$ solver.
Although, these two solves can be performed either in an adiditive or a
multiplicative manner, we found that the reduction of pressure iterations
resulting from the latter is essential for the new solver to be competitive
with state of the art AMG solvers.

There are multiple aspects in the current solver that can be improved in future
research work to make it more robust and faster.
Currently, we are solving both local and global coarse problems of the Schwarz
method in double precision.
Switching to single precision will cut the communication cost by a factor of
two thus could result in a considerable speed up in global coarse solve time.
Similarly, local Cholesky solve time can see a speedup too due to the low
memory footprint of the Cholesky factor.
Another effective way to reduce memory footprint and communication cost is to
reduce nonzero entries in both the local Cholesky factor as well as the
$XX^T$ factorization with a reordering of the respective sparse systems.
Furthermore, reduction of the nonzeros in the factors reduce the number
of floating point operations that has to be performed during a solve.
Another avenue worth exploring is the use of inexact solvers for both local
Schwarz problem as well as the global coarse grid problem.
This could be a winning strategy if the increase of pressure iterations and
subsequent increase of the time spent in finer levels of $p$-multigrid can be
compensated by the reduction in time spent in the coarse solve.
For the local problem, incomplete Cholesky on the CPU or an iterative method
like CG with a suitable preconditioner on GPU can be used as an inexact solver.


\section*{Acknowledgments}

This material is based upon work supported by the U.S. Department of Energy,
Office of Science, under contract DE-AC02-06CH11357 and by the Exascale Computing
Project (17-SC-20-SC). The research used resources at the Oak Ridge Leadership
Computing Facility at Oak Ridge National Laboratory, which is supported by
the Office of Science of the U.S. Department of Energy under Contract DE-AC05-00OR22725.




\bibliographystyle{cas-model2-names}
\addtolength{\baselineskip}{-.1\baselineskip}
\bibliography{bibs/emmd}
\addtolength{\baselineskip}{+.111111111\baselineskip}

\bio{}
\endbio

\bio{}
\endbio

\end{document}